\documentclass{article}
\usepackage{amsmath,amssymb}

\newtheorem{theorem}{Theorem}[section]
\newtheorem{lemma}[theorem]{Lemma}
\newtheorem{proposition}[theorem]{Proposition}
\newtheorem{corollary}[theorem]{Corollary}

\newcommand{\ch}{\operatorname{char}}

\newcommand{\lin}{\operatorname{lin}}

\newcommand{\supp}{\operatorname{supp}}
\newcommand{\gr}{\operatorname{gr}}

\newcommand{\Alt}{\operatorname{Alt}}
\newcommand{\Sym}{\operatorname{Sym}}
\newcommand{\FM}{\operatorname{FM}}
\newcommand{\FaM}{\operatorname{FaM}}

\newenvironment{proof}{\par\noindent{\bf Proof.}}{$\qed$\par\bigskip}
\newcommand{\qed}{\enspace\vrule  height6pt  width4pt  depth2pt}
\begin{document}

\title{Algebras and groups defined by permutation relations of
alternating type \thanks{ Research  partially
supported by grants of MICIN-FEDER (Spain)
MTM2008-06201-C02-01, Generalitat de Catalunya
2005SGR00206, Onderzoeksraad of Vrije Universiteit
Brussel, Fonds voor Wetenschappelijk Onderzoek
(Belgium), Flemish-Polish bilateral agreement
BIL2005/VUB/06 and MNiSW research grant N201 004
32/0088 (Poland).}}

\author{Ferran Ced\'o \and Eric Jespers \and Jan Okni\'nski}
\date{}
\maketitle

\begin{abstract}
The class of finitely presented algebras over a
field $K$ with a set of generators $a_{1},\ldots
, a_{n}$ and defined by homogeneous relations of
the form
 $a_{1}a_{2}\cdots a_{n} =a_{\sigma (1)} a_{\sigma (2)} \cdots
 a_{\sigma (n)}$,
where $\sigma$ runs through $\Alt_{n}$, the
alternating group, is considered. The associated
group, defined by the same (group) presentation,
is described.  A description of the radical of
the algebra is found. It turns out that the
radical is a finitely generated ideal that is
nilpotent and it is determined by a congruence on
the underlying monoid, defined by the same
presentation.
\end{abstract}


\section{Introduction}\
In recent literature a lot of attention is given
to concrete classes of finitely presented
algebras $A$ over a field $K$ defined by
homogeneous semigroup relations, that is,
relations of the form $w=v$, where $w$ and $v$
are words of the same length in a generating set
of the algebra. Of course such an algebra is a
semigroup algebra $K[S]$, where $S$ is the monoid
generated by the same presentation. Particular
classes show up in different areas of research.
For example, algebras yielding set theoretic
solutions of the Yang-Baxter equation (see for
example
\cite{eting,gat-sol,gat-van,jes-okn-Itype,rump})
or algebras related to Young diagrams,
representation theory and algebraic combinatorics
(see for example
\cite{chinese,plactic,fulton,jaszu-okn,las-lec}).
In all the mentioned algebras there are  strong
connections between the structure of the algebra
$K[S]$, the underlying semigroup $S$ and the
underlying group $G$, defined by the same
presentation as the algebra.

In \cite{alghomshort} the authors introduced and initiated a study
of combinatorial and algebraic aspects of the following  new class
of finitely presented algebras over a field $K$: $$A= K\langle
a_{1},a_{2},\ldots , a_{n} \mid a_{1}a_{2}\cdots a_{n} =a_{\sigma
(1)} a_{\sigma (2)} \cdots
 a_{\sigma (n)},\; \sigma \in H\rangle ,$$
where $H$ is a subset of the symmetric group
$\Sym_{n}$ of degree $n$. So $A=K[S_{n}(H)]$
where $$S_{n}(H)=\langle a_{1},a_{2},\ldots ,
a_{n} \mid a_{1}a_{2}\cdots a_{n} =a_{\sigma (1)}
a_{\sigma (2)} \cdots
 a_{\sigma (n)},\; \sigma \in H\rangle ,$$
the monoid with the ``same'' presentation as the
algebra. By $G_{n}(H)$ we denote the group
defined by this presentation. So $$G_{n}(H)=\gr (
a_{1},a_{2},\ldots , a_{n} \mid a_{1}a_{2}\cdots
a_{n} =a_{\sigma (1)} a_{\sigma (2)} \cdots
 a_{\sigma (n)},\; \sigma \in H) .$$
Two obvious examples are: the free $K$-algebra
$K[S_{n}(\{ 1\})]=K\langle a_{1},\ldots ,
a_{n}\rangle $ with $H=\{1\}$ and $S_{n}(\{ 1
\})=\FM_{n}$, the rank $n$ free monoid, and the
commutative polynomial algebra
$K[S_{2}(\Sym_{2})]=K[a_{1},a_{2}]$ with
$H=\Sym_{2}$  and $S_{n}(H)=\FaM_{2}$, the rank
$2$ free abelian monoid. For $M=S_{n}(\Sym_{n})$,
the latter can be extended as follows
(\cite[Proposition~3.1]{alghomshort}): the
algebra $K[M]$ is the subdirect product of the
commutative polynomial algebra $K[a_{1},\ldots
,a_{n}]$ and a primitive monomial algebra that is
isomorphic to $K[M]/K[Mz]$, with
$z=a_{1}a_{2}\cdots a_{n}$, a central element.

On the other hand, let $M=S_{n}(H)$ where $H=\gr (\{ (1\, ,2\, ,
\ldots ,n)\} )$, a cyclic group of order $n$. Then
(\cite[Theorem~2.2]{alghomshort}) the monoid $M$ is cancellative
and it has a group $G$ of fractions of the form $G=M\langle
a_{1}\cdots a_{n}\rangle ^{-1}\cong F\times C$, where $F=\gr
(a_{1},\ldots , a_{n-1})$ is a free group of rank $n-1$ and $C=\gr
(a_{1}\cdots a_{n})$ is a cyclic infinite group. The algebra
$K[M]$ is a domain and it is semiprimitive. Moreover
(\cite[Theorem~2.1]{alghomshort}), a normal form of elements of
the algebra can be given. It is worthwhile mentioning that the
group $G$ is an example of a cyclically presented group. Such
groups arise in a very natural way as fundamental groups of
certain $3$-manifolds \cite{edjvet}, and their algebraic structure
also receives a lot of attention; for a recent work and some
references see for example \cite{cavic}.

In this paper we continue the investigations on the algebras
$K[S_{n}(H)]$ and the groups $G_{n}(H)$. First we will prove some
general results and next will give a detailed account in case $H$
is the alternating group $\Alt_{n}$ of degree $n$. It turns out
that the structure of the group $G_{n}(H)$ can be completely
determined and the algebra $K[S_{n}(H)]$ has some remarkable
properties. In order to state our main result we fix some
notation. Throughout the paper $K$ is a field. If $b_{1},\ldots,
b_{m}$ are elements of a monoid $M$ then we denote by $\langle
b_{1},\ldots , b_{m}\rangle$ the submonoid generated by
$b_{1},\ldots , b_{m}$. If $M$ is a group then $\gr (b_{1},\ldots
,b_{m})$ denotes the subgroup of $M$ generated by $b_{1},\ldots,
b_{m}$. Clearly, the defining relations of an arbitrary $S_{n}(H)$
are homogeneous. Hence, it has a natural degree or length
function. This will be used freely throughout the paper. By $\rho
= \rho_{S}$ we denote the least cancellative congruence on a
semigroup $S$.  If $\eta $ is a congruence on $S$ then
 $I(\eta)= \lin \{ s-t \mid s,t\in M, (s,t)\in \eta\}$ is
 the kernel of the natural epimorphism
$K[S]\longrightarrow K[S/\eta ]$. For a ring $R$, we denote by
$\mathcal{J}(R)$ its Jacobson radical and by $\mathcal{B}(R)$ its
prime radical. Our main result reads as follows.

\begin{theorem} \label{main-alt}
Suppose $K$ is a field and $n\geq 4$. Let
$M=S_{n}(\Alt_{n})$, $z=a_1a_2\cdots a_n\in M$
and
 $G=G_{n}(\Alt_{n})$. The following properties
 hold.
\begin{itemize}
\item[(i)] $C=\{1,a_{1}a_{2}a_{1}^{-1}a_{2}^{-1}\}$ is a nontrivial
central subgroup of $G$ and $G/C$ is a free
abelian group of rank $n$. Moreover $D=\gr
(a_{i}^{2} \mid i=1,\ldots ,n)$ is a central
subgroup of $G$ with $G/(CD)\cong
(\mathbb{Z}/2\mathbb{Z})^{n}$.
\item[(ii)] $K[G]$ is a noetherian algebra satisfying a polynomial
identity (PI, for short). If $K$ has
characteristic $\neq 2$, then
$\mathcal{J}(K[G])=0$. If $K$ has characteristic
$2$, then
$\mathcal{J}(K[G])=(1-a_1a_2a_1^{-1}a_2^{-1})K[G]$
and $\mathcal{J}(K[G])^2= 0$.
\item[(iii)] The element $z^2$ is central in $M$ and $z^{2}M$ is a
cancellative ideal of $M$ such that $G \cong (z^{2}M)\langle
z^{2}\rangle ^{-1}$. Furthermore, $K[M/\rho]$ is a noetherian
PI-algebra and $\mathcal{J}(K[M])$ is nilpotent.
\item[(iv)] Suppose $n$ is odd. Then $z$ is central in $M$ and
$0\neq \mathcal{J}(K[M])=I(\eta)$ for a
congruence $\eta $ on $M$  and
$\mathcal{J}(K[M])$ is a finitely generated
ideal.
\item[(v)]  Suppose $n$ is even and $n\geq 6$. If
 $K$ has characteristic $\neq 2$,
then $\mathcal{J}(K[M])=0$. If $K$ has
characteristic $2$, then $0\neq
\mathcal{J}(K[M])=I(\eta)$ for a congruence
$\eta$ on $M$ and $\mathcal{J}(K[M])$ is a
finitely generated ideal.
\end{itemize}
\end{theorem}

Part $(v)$ of Theorem \ref{main-alt} is also true for $n=4$, but
its proof is quite long for this case and it requires additional
technical lemmas. (The interested reader can find a proof of this
in \cite{alt4algebra}).

So, in particular,  the Jacobson radical is determined by a
congruence relation on the semigroup $S_{n}(\Alt_{n})$, it is
nilpotent and finitely generated as an ideal. In
\cite{alghomshort} the question was asked whether these properties
hold for all algebras $K[S_{n}(H)]$, for subgroups $H$ of
$\Sym_{n}$.

\section{General results}

In this section we prove some preparatory general
properties of the monoid algebra $K[S_n(H)]$ for
an arbitrary subset $H$ of $\Sym_{n}$ with $n\geq
3$. To simplify notation, throughout this section
we put
 \begin{eqnarray}
 M&=&\langle a_1,a_2,\dots ,
a_n \mid a_1a_2\cdots
a_n=a_{\sigma(1)}a_{\sigma(2)}\cdots
a_{\sigma(n)}, \;  \sigma\in H\rangle
\label{monoidpresentation}
\end{eqnarray}
If $\alpha =\sum_{x\in M} k_{x}x\in K[M]$, with
each $k_{x}\in K$, then the finite set $\{ x\in
M\mid k_{x}\neq 0\}$ we denote by $\supp (\alpha
)$. It is called the support of $\alpha$.

\begin{proposition}\label{transitive}
Suppose that there exists $k$ such that, for all
$\sigma\in H$, $\sigma(1)\neq k$ and
$\sigma(n)\neq k$. Then $\mathcal{J}(K[M])=0$.
\end{proposition}

\begin{proof}
Suppose that $\mathcal{J}(K[M])\neq 0$. Let
$\alpha\in \mathcal{J}(K[M])$ be a nonzero
element. Because, by assumption, $\sigma (1)\neq
k$ for all $\sigma \in H$, we clearly get that
$a_k^2\alpha\neq 0$. As $a_k^2\alpha \in
\mathcal{J}(K[M])$, there exists $\beta\in K[M]$
such that $a_k^2\alpha+\beta+\beta
a_k^2\alpha=0$. Obviously, $\beta\notin K$. Let
$\alpha_1$, $\beta_1$ be the homogeneous
components (for the natural
$\mathbb{Z}$-gradation of $K[M]$) of $\alpha$ and
$\beta$ of maximum degree respectively. Then
$\beta_1a_k^2\alpha_1=0$. In particular, there
exist $w_1,w_2$ in the support of $\beta_1$ and
$w_1',w_2'$ in the support of $\alpha_1$ such
that $$w_1a_k^2w_1'=w_{2}a_k^2w_2'$$ and either
$w_1\neq w_2$ or $a_k^2w_1'\neq a_k^2w_2'$. But,
because $\sigma (n)\neq k$ for all $\sigma \in
H$, this is impossible. Therefore
$\mathcal{J}(K[M])=0$.
\end{proof}

\begin{corollary}
If $H$ is a subgroup of $\Sym_n$ and
$\mathcal{J}(K[M])\neq 0$ then $H$ is a
transitive subgroup of $\Sym_{n}$.
\end{corollary}

\begin{proof}
Suppose that $H$ is a subgroup of $\Sym_n$ and
$\mathcal{J}(K[M])\neq 0$. By
Proposition~\ref{transitive}, for all $k$ there
exist $\sigma\in H$ such that either
$\sigma(1)=k$ or $\sigma(n)=k$. Suppose that $H$
is not transitive. Then there exists $1\leq j\leq
n$ such that $j\notin\{ \sigma(1)\mid \sigma\in H
\}$. Hence there exists $\sigma\in H$ such that
$\sigma (n)=j$. Thus the orbits $I_1=\{
\sigma(1)\mid \sigma\in H \}$ and $I_2=\{
\sigma(n)\mid \sigma\in H \}$ are disjoint
nonempty sets such that $I_1\cup I_2=\{ 1,2,\dots
,n\}$.  So, there are no defining relations of
the form $a_{1} \cdots =a_{n}\cdots $, nor of the
form $\cdots a_{1} = \cdots a_{n}$. Consequently,
if $0\neq \alpha \in K[M]$ then $a_{n}^{2}\alpha
\neq 0$ and $\alpha a_{1}^{2}\neq 0$.

Let $\alpha\in \mathcal{J}(K[M])$ be a nonzero
element. Then, $a_1^2a_n^2\alpha\neq 0$, and
there exists $\beta\in K[M]$ such that
$a_1^2a_n^2\alpha+\beta+\beta
a_1^2a_n^2\alpha=0$. Clearly, it follows that
$\beta\notin K$. Let $\alpha_1$, $\beta_1$ be the
homogeneous components of $\alpha$ and $\beta$ of
maximum degree respectively. We obtain that
$\beta_1a_1^2a_n^2\alpha_1=0$. In particular,
there exist $w_1,w_2$ in the support of $\beta_1$
and $w_1',w_2'$ in the support of $\alpha_1$ such
that $(w_{1},w_{1}')\neq (w_{2},w_{2}')$ and
$$w_1a_1^2a_n^2w_1'=w_2a_1^2a_n^2w_2'.$$ Again,
because there are no defining relations of the
form $a_{1} \cdots =a_{n}\cdots $ nor of the form
$\cdots a_{1} = \cdots a_{n}$, this yields a
contradiction. Therefore $H$ is transitive.
\end{proof}

 Let $z=a_1a_2\cdots a_n\in M$.  The fact that $z$ is central
in $M=S_{n}(H)$ for the case of the cyclic group
$H$ generated by $(1,2,\ldots ,n)$ was an
important tool in \cite{alghomshort}. In
Section~\ref{altern} we will show that $z^{2}$ is
central if $M=S_{n}(\Alt_{n})$. In the following
properties we show that the centrality of
$z^{m}$, for some positive integer $m$, has some
impact on the algebraic structure of $M$ and
$K[M]$. But first we determine when $z$ is
central in case $H$ is a subgroup of $\Sym_n$.

\begin{proposition} \label{centralityofz} Suppose $H$ is a subgroup of
$\Sym_{n}$ and put $z=a_{1}a_{2}\cdots a_{n}$.
The following conditions are equivalent.
\begin{itemize}
\item[(i)] $z$ is central in $M=S_{n}(H)$,
\item[(ii)] $a_{1}z=za_{1}$,
\item[(iii)] $H$ contains the subgroup of
$\Sym_n$ generated by the cycle $(1,2,\dots ,n)$.
\end{itemize}
\end{proposition}
\begin{proof} Let $H_0$ denote the subgroup of $\Sym_n$ generated
by the cycle $(1,2,\dots , n)$. Assume $H_{0}\subseteq H$. Then
$M=S_{n}(H)$ is an epimorphic image of $S_{n}(H_0)$.  As
$a_{1}a_{2}\cdots a_{n}$ is central in $S_{n}(H_0)$, it follows
that $z$ indeed is central in $M$.

Assume now that $a_{1}z=za_{1}$. We need to show
that $H_{0}\subseteq H$. Every defining relation
can be written in the form: $z=a_{k}c_{k}$, with
$1\leq k \leq n$, $c_{k}=\prod_{i=1,\; i\neq k}
a_{\tau (i)}$, $\tau \in\Sym (\{ 1,\dots , n\}
\setminus \{ k \}) \subseteq \Sym_{n}$. By
assumption, $a_{1}^{2} a_{2}\cdots
a_{n}=a_{1}z=za_{1}=a_{k}c_{k}a_{1}$. Since
$a_{k}c_{k}$ is a product of distinct generators,
there must exist a relation of the form
$c_{k}a_{1}=z$. Since also $z$ is a product of
distinct generators, it follows that $k=1$. Thus
$z=a_{1}c_{1}=c_{1}a_{1}$. The former equality
yields that $\sigma_{1}
=\left(\begin{array}{ccccc} 1&2&\dots&n-1&n\\
1&\tau(2)&\dots&\tau(n-1)&\tau(n)
\end{array}\right)\in H$ and the equality
$z=c_{1}a_{1}$ gives $\sigma_{2} = \left(
\begin{array}{ccccc}
1&2&\dots&n-1&n\\ \tau(2)&\tau(3)&\dots&\tau(n)&1
\end{array}\right)\in H$. Hence $(1,2,\dots ,
n)=\sigma_{1}\sigma_{2}^{-1}\in H$ and so
$H_{0}\subseteq H$. The result follows.
\end{proof}

Assume now that $z^{m}$ is central, for some positive integer $m$.
Note that then  the binary relation $\rho'$ on $M$, defined by
$s\rho' t$ if and only if there exists a nonnegative integer $i$
such that $sz^{i}=tz^{i}$, is a congruence on $M$.  We now show
that $G_{n}(H)$ is the group of fractions of
$\overline{M}=M/\rho'$. We denote by $\overline{a}$ the image in
$\overline{M}$ of $a\in M$ under the natural map $M\longrightarrow
\overline{M}$.

\begin{lemma}  \label{canc}  Suppose that $z^m$
is central for some positive integer $m$. Then, $\rho'=\rho$ is
the least cancellative congruence on $M$ and $Ma \cap aM \cap
\langle z^{m} \rangle \neq \emptyset$ for every $a\in M$.

In particular, $\overline{M}=M/\rho$ is a
cancellative monoid and $G=\overline{M}\langle
\overline{z}^{-m}\rangle$ is the group of
fractions of $\overline{M}$. Moreover, $G\cong
G_{n}(H)=\gr (a_{1},\ldots , a_{n} \mid
a_{1}a_{2}\cdots a_{n} =a_{\sigma (1)}a_{\sigma
(2)} \cdots a_{\sigma (n)},\; \sigma \in H)$.
\end{lemma}
\begin{proof}
Since $z^m$ is central, we already know that the binary relation
$\rho'$ is a congruence on $M$. Let $a=a_{i_1}a_{i_2}\cdots
a_{i_k}\in M$. We shall prove that $aM\cap\langle z^m\rangle\neq
\emptyset$ by induction on $k$. For $k=0$, this is clear. Suppose
that $k>0$ and that $bM\cap \langle z^m\rangle\neq \emptyset$ for
all $b\in M$ of degree less than $k$. Thus there exists $r\in M$
such that $a_{i_1}\cdots a_{i_{k-1}}r\in \langle z^m\rangle$.
Since $a_{i_{k}}z^{m}=z^{m}a_{i_{k}}$, it follows easily from the
type of the defining relations for $M$ that there exists $w\in M$
such that $a_{i_k}w=z$. We thus get that $awz^{m-1}r=a_{i_1}\dots
a_{i_{k-1}}z^m r=a_{i_1}\dots a_{i_{k-1}}rz^m\in \langle
z^m\rangle$. Similarly we see that $Ma\cap\langle z^m \rangle\neq
\emptyset$. Therefore $\rho'$ is the least cancellative congruence
on $M$ and $\overline{M}\langle \overline{z}^{-m}\rangle$ is the
group of fractions of $\overline{M}$ and the second assertion also
follows.
\end{proof}

\begin{proposition}\label{locnil}  Suppose that $z^{m}$
is central for some positive integer $m$.  Let
$\alpha_1,\ldots ,\alpha_k\in I(\rho ) \cap
K[Mz^{m}]$. Then the ideal $\sum_{i=1}^{k}
K[M]\alpha_i K[M]$ is nilpotent. In particular,
$I(\rho )\cap K[Mz^{m}] \subseteq
{\mathcal{B}}(K[M])$.
\end{proposition}
\begin{proof}
  Let $\alpha_1,\dots,\alpha_k\in I(\rho ) \cap
K[Mz^{m}]$. Clearly there exists a positive
integer $N$ such that $\alpha_i z^{mN}=0$, for
all $i=1,\dots ,k$. Since $z^{m}$ is central and
$\alpha_i\in K[Mz^{m}]$, we have that
$(\sum_{i=1}^mK[M]\alpha_iK[M])^{N+1}=0$  and the
result follows.
\end{proof}

\begin{proposition} \label{radicalgeneral}
The
following properties hold.
\begin{itemize}
\item[(i)] $\mathcal{J}(K[M]/K[MzM])=0$.
\item[(ii)] $\mathcal{J}(K[M])\subseteq K[Mz \cup zM]$.
\item[(iii)] $\mathcal{J}(K[M])^{3} \subseteq K[Mz^{2}MzM \cup MzMz^{2}M] \subseteq K[Mz^{2}M]$.
\end{itemize}
If, furthermore, $z^{m}$ is central for some
positive integer $m$, $Mz^{k}M$ is cancellative
and $\ch (K)=0$ then $K[Mz^kM]$ has no nonzero
nil ideal. In particular,
$\mathcal{B}(K[Mz^{k}M])=0$. Furthermore,  if
$k=2$ then $\mathcal{B}(K[M])^{3}=0$.
\end{proposition}
\begin{proof}
To prove the first part, let $X$ be the free monoid
with basis $x_1, x_2,\dots, x_n$. Then
$$K[M]/K[MzM]\cong K[X]/K[J],$$ where
$J=\bigcup_{\sigma\in H\cup\{ 1\}}
Xx_{\sigma(1)}\cdots x_{\sigma(n)}X$.  Note that
$X/J$ has no nonzero nilideal. Hence, by
\cite[Corollary~24.7]{book}, $K[M]/K[MzM]$ is
semiprimitive. Therefore
$\mathcal{J}(K[M]/K[MzM])=0$.

To prove the second and third part,
suppose that $\alpha
=\sum_{i=1}^{m}\lambda_{i}s_{i}\in
\mathcal{J}(K[M])$, with $\supp (\alpha ) =\{
s_{1},\ldots , s_{m}\}$ of cardinality $m$ and
$\lambda_{i}\in K$,  is a homogeneous element
(with respect to the gradation defined by the
natural length function on $M$). Then $\alpha$ is
nilpotent (see for example
\cite[Theorem~22.6]{passmancrossed}). Suppose
that $s_{1}\not \in zM$ and $s_{1}\not \in Mz$.
Let $i,j$ be such that $s_{1}\in a_{i}M\cap
Ma_{j}$. Then, for every $k\geq 1$, the element
$(s_{1}a_{j}a_{i})^{k}=s_{1}a_{j}a_{i}s_{1}a_{j}a_{i}
\cdots $ can only be rewritten in $M$ in the form
$(s'a_{j}a_{i})^{k}$, where $s'\in M$ is such
that $s'=s_{1}$. Therefore, $\alpha a_{j}a_{i}\in
\mathcal{J}(K[M])$ is not nilpotent, a
contradiction. It follows that $s_{1}$, and
similarly every $s_{i}\in Mz\cup zM$. Again by
\cite[Theorem~22.6]{passmancrossed}), we know
that $\mathcal{J}(K[M])$ is a homogeneous ideal.
This implies that $\mathcal{J}(K[M])\subseteq
K[Mz\cup zM]$. Hence
$\mathcal{J}(K[M])^{3}\subseteq
Mz\mathcal{J}(K[M])zM \subseteq K[Mz^{2}MzM \cup
MzMz^{2}M]$. This finishes the proof of
statements $(ii)$ and $(iii)$.

To prove the last part, assume $\ch (K)=0$,
$Mz^{k}M$ is cancellative and $z^{m}$ is central
for some positive integer $m$. Since
$a_{i}z^{m}=z^{m}a_{i}$, it follows from the type
of the  defining relations for $M$ that $z\in
a_{i}M\cap Ma_{i}$ for every $1\leq i \leq n$.
Hence, by Lemma~\ref{canc}, we know that
$Mz^{k}M$ has a group of fractions $G$ (that is
obtained by inverting the powers of  the central
element $z^{km}$).  Let $I$ be a nil ideal of
$K[Mz^kM]$. Then $K[G]IK[G]=I\langle
z^{-km}\rangle$ is a nil ideal of $K[G]$. Since,
by assumption, $\ch (K)=0$, we know from
\cite[Theorem~2.3.1]{passman} that then $I=0$.
So, if $k=2$ then, by the first part of the
result, $\mathcal{B}(K[M])^{3}\subseteq
K[Mz^{2}M] \cap \mathcal{B}(K[M])$.  Since
$K[Mz^{2}M] \cap \mathcal{B}(K[M])$ is a nil
ideal of $K[Mz^2M]$, the result follows.
\end{proof}

\begin{corollary}\label{radical}
Suppose $z$ is central. The following properties
hold.
\begin{itemize}
\item[(i)] If $\mathcal{J}(K[\overline{M}])=0$ then $\mathcal{J}(K[M])=I(\rho )\cap K[Mz]$.
\item[(ii)] If $\mathcal{B}(K[\overline{M}])=0$ then
$\mathcal{B}(K[M])=I(\rho )\cap K[Mz]$.
\item[(iii)] If $\mathcal{B}(K[M])= 0$ then $Mz$ is cancellative.
The converse holds provided $\ch (K)=0$.
\end{itemize}
\end{corollary}

\begin{proof}
$(i)$  By Proposition~\ref{radicalgeneral},
$\mathcal{J}(K[M])\subseteq K[Mz]$. Note that
$K[\overline{M}]=K[M/\rho]=K[M]/I(\rho )$. Hence,
if $\mathcal{J}(K[\overline{M}])=0$, we get that
$\mathcal{J}(K[M])\subseteq I(\rho ) \cap K[Mz]$.
By Proposition~\ref{locnil} we thus obtain that
$\mathcal{J}(K[M])=I(\rho )\cap K[Mz]$.

$(ii)$ If $K[\overline{M}]$ is semiprime, then,
by  Proposition~\ref{radicalgeneral},
$\mathcal{B}(K[M])\subseteq I(\rho )\cap K[Mz]$.
Thus, by Proposition~\ref{locnil},
$\mathcal{B}(K[M])=I(\rho )\cap K[Mz]$.

$(iii)$ Because of Proposition~\ref{locnil},
we know that $I(\rho ) \cap K[Mz] \subseteq
\mathcal{B}(K[M])$. Suppose now that
$\mathcal{B}(K[M])= 0$. Then, $\rho$ restricted
to $Mz$ must be the trivial relation, i.e., $Mz$
is cancellative. Conversely, assume that $\ch
(K)=0$ and $Mz$ is cancellative. Then, by
Proposition~\ref{radicalgeneral},
$\mathcal{B}(K[M])$ is a nil ideal of $K[Mz]$,
and thus $\mathcal{B}(K[M])=0$, as desired.
\end{proof}

\section{The monoid $S_{n}(\Alt_{n})$}

In this section we investigate the monoid
$S_{n}(\Alt_{n})$
with $n\geq 4$.  The information obtained is
essential to prove our  main result,
Theorem~\ref{main-alt}.  Note that the cycle $(1,
2, \ldots , n)\in \Alt_{n}$ if and only if $n$ is
odd. Hence by Proposition~\ref{centralityofz},
$z=a_{1}a_{2}\cdots a_{n}$ is central if and only
if $n$ is odd.
However, for arbitrary $n$, we will show that
$z^{2}$ is central and that the ideal
$S_{n}(\Alt_{n})z^{2}$ is cancellative as a
semigroup and we also will determine the
structure of its group of fractions
$G_{n}(\Alt_{n})$. This information will be
useful to determine the radical of the algebra
$K[S_{n}(\Alt_{n})]$.

Throughout this section $n\geq 4$,
$M=S_n(\Alt_n)$ and $G=G_n(\Alt_n)$. Let
$\sigma\in \Alt_n$. Since the set of defining
relations of $M$ (of $G$, respectively) is
$\sigma$-invariant, $\sigma$ determines the
automorphism of $M$ (of $G$ respectively) defined
by $\sigma (a_{i_1}^{n_1}\cdots
a_{i_m}^{n_m})=a_{\sigma(i_1)}^{n_1}\cdots
a_{\sigma(i_m)}^{n_m}$.

\begin{lemma}\label{ij}
Let $z=a_1a_2\cdots a_n\in M$.
\begin{itemize}
\item[(i)] If $n\geq 4$ then $a_ia_jz=za_ia_j$, for any different integers $1\leq i,j\leq n$.
\item[(ii)] If $n\geq 5$ then $a_ia_ja_kz=a_ja_ka_iz$ and $za_ia_ja_k=za_ja_ka_i$, for
any three different integers $1\leq i,j,k\leq n$.
\item[(iii)] If $n=4$ and $1\leq i,j,k\leq n$ are
three different integers then
\begin{enumerate}
 \item if $a_{i}a_{j}a_{k}a_{l}=z$ then $a_ia_ja_kz=a_ja_ka_iz=a_ka_ia_jz=za_ka_ja_i$,
 \item if $a_{l}a_{i}a_{j}a_{k}=z$ then $za_ia_ja_k=za_ja_ka_i=za_ka_ia_j=a_{k}a_ja_iz$.
\end{enumerate}
\item[(iv)] If $n\geq 6$ is even then $a_ia_ja_kz=za_ja_ia_k$, for
any three different integers $1\leq i,j,k\leq n$.
\end{itemize}
\end{lemma}
\begin{proof}
$(i)$ If $1\leq i, j\leq n$ are different then there
exists $\sigma\in \Alt_n$ such that $\sigma(1)=i$
and $\sigma(2)=j$. Hence
\begin{eqnarray*}
a_ia_jz&=&a_ia_j\sigma(1,2,\dots
,n)^2(a_1a_2\cdots a_n)\\
&=&a_ia_ja_{\sigma(3)}\cdots
a_{\sigma(n)}a_{\sigma(1)}a_{\sigma(2)}=za_ia_j.
\end{eqnarray*}

$(ii)$ and $(iv)$  Suppose that $n\geq 5$. In
this case, for any three different integers
$1\leq i,j,k\leq n$ there exists $\sigma \in
\Alt_{n}$ such that $\sigma (1)=i, \sigma (2)=j,
\sigma (3)=k$. Let $\tau=\tau_n\in \Sym_n$ be
defined by $\tau =id$ if $n$ is odd, and $\tau =
(i,j)$ if $n$ is even. So $\tau \sigma
(1,2,\ldots  , n)^{3}\in \Alt_{n}$.
Hence in $M$ we get
\begin{eqnarray*}
a_{i}a_{j}a_{k}z & = & a_ia_ja_k\tau\sigma
(1,2,\dots ,n)^3 (a_{1}a_{2}\cdots a_{n})\\ & = &
(a_{\sigma(1)}a_{\sigma(2)}a_{\sigma(3)})(a_{\sigma(4)}\cdots
a_{\sigma(n)}a_{\tau(i)}a_{\tau(j)}a_k)
\\
 &=& \sigma (z) a_{\tau (i)} a_{\tau (j)} a_{k}.
\end{eqnarray*}
 In particular, $(iv)$ follows. Since
$(1,2,3)\in \Alt_{n}$, this yields
\begin{eqnarray*}
 a_{i}a_{j}a_{k}z & = & (\sigma(1,2,3)(a_1a_2\cdots
a_n))a_{\tau(i)}a_{\tau(j)}a_k
\\
& = & (a_{j}a_{k}a_{i})a_{\sigma(4)}\cdots
a_{\sigma(n)}a_{\tau(i)}a_{\tau(j)}a_k=a_{j}a_{k}a_{i}z
.
\end{eqnarray*}
Similarly one proves that
\begin{eqnarray*} za_ia_ja_k=za_ja_ka_i,
\end{eqnarray*}
for $n\geq 5$.

$(iii)$ Suppose that $n=4$. Let $\{ i,j,k,l\}=\{
1,2,3,4\}$. Then either $a_ia_ja_ka_l=z$ or
$a_la_ia_ja_k=z$. If $a_ia_ja_ka_l=z$, then
$$z=a_ia_ja_ka_l=a_ja_ka_ia_l=a_ka_ia_ja_l,$$ and,
since $z\in a_lM$, we get
$$a_ia_ja_kz=a_ja_ka_iz=a_ka_ia_jz.$$ Clearly,
$a_ia_ja_kz= a_{i}a_{j}a_{k}(a_{l}a_{k}a_{j}a_{i})=
za_{k}a_{j}a_{i}$.

 Similarly,
if $a_la_ia_ja_k=z$, we get
$$za_ia_ja_k=za_ja_ka_i=za_ka_ia_j=(a_{k}a_{j}a_{i}a_{l})a_ka_ia_j=a_{k}a_{j}a_{i}z.$$
\end{proof}

\begin{lemma}\label{z2central}
Let $z=a_1a_2\cdots a_n\in M$. Then $z^2$ is
central in $M$.
\end{lemma}

\begin{proof}
If  $n\geq 6$ and $n$ is even then
\begin{eqnarray*}
z^2a_1&=&za_1a_2\cdots
a_na_1=a_1a_2a_3a_4za_5\cdots a_na_1\quad
(\mbox{by Lemma~\ref{ij}.(i)})\\
&=&a_1a_2a_3a_4((1,5)(2,3)(2,3,\dots ,
n)^3(a_1a_2\cdots a_n))a_5\cdots a_na_1\\
&=&a_1a_2a_3a_4(a_5a_1a_6\cdots
a_na_3a_2a_4)a_5\cdots a_na_1\\
&=&a_1((1,2,3,4,5)(a_1a_2\cdots
a_n))(2,3)(1,2,\dots ,n)(a_1a_2\dots a_n)\\
&=&a_1z^2.
\end{eqnarray*}
If  $n=4$ then
\begin{eqnarray*}
a_1 z^2 &=& a_1(a_3a_4a_1a_2)z=a_1a_3a_4z a_1a_2
\quad\mbox{(by Lemma \ref{ij})}\\ &=&
a_1a_3a_4(a_2a_3 a_1a_4) a_1a_2 = za_3
a_1a_4a_1a_2\\ &=& a_3a_1za_4a_1a_2\quad\mbox{(by
Lemma \ref{ij})}\\ &=&
a_3a_1(a_2a_1a_4a_3)a_4a_1a_2=a_3a_1a_2a_1a_4(a_3a_2a_4a_1)\\
&=&a_3a_1za_2a_4a_1=za_3a_1a_2a_4a_1\quad\mbox{(by
Lemma \ref{ij})}\\ &=&z^2a_1.
\end{eqnarray*}
Since $\Alt_{n}$ is transitive, we get that
$z^{2}$ is central for all even $n$.
Since $z$ is central in $M$ for all odd $n$, the
assertion follows.
\end{proof}

\begin{lemma}\label{central4}
For $n=4$, $a_1a_2a_4a_3z=\sigma(a_1a_2a_4a_3)z$,
for all $\sigma\in \Alt_4$, and it is central in
$M$. In particular, $\sigma (z) z =z\sigma (z) =
z \gamma (z)=\gamma (z)  z$ for any $\sigma
,\gamma  \in \Sym_{4}$ of the same parity.
\end{lemma}

\begin{proof}
By Lemma~\ref{ij}, we have
$$a_1a_2a_4a_3z=a_1(a_3a_2a_4)z=a_1(a_4a_3a_2)z,$$
and also
\begin{eqnarray*}
a_1a_2a_4a_3z&=&za_1a_2a_4a_3=z(a_2a_4a_1)a_3=a_2a_4a_1a_3z,\\
a_1a_2a_4a_3z&=&za_1a_2a_4a_3=z(a_4a_1a_2)a_3=a_4a_1a_2a_3z,\\
a_1a_2a_4a_3z&=&a_1a_3a_2a_4z=za_1a_3a_2a_4=z(a_3a_2a_1)a_4=a_3a_2a_1a_4z.
\end{eqnarray*}
Thus $a_1a_2a_4a_3z=\sigma(a_1a_2a_4a_3)z$ for
all $\sigma\in \Alt_4$. In particular, $\sigma
(z) z =\gamma (z) z$ for odd permutations $\sigma
, \gamma$. Of course such an equality also holds
if $\gamma ,\sigma$ are even. Note that, because
of Lemma~\ref{ij}, $z\sigma (z) =\sigma (z) z$
for any permutations $\sigma$.

In order to prove that $a_1a_2a_4a_3z$ is central
we only need to show that $a_1a_2a_4a_3z a_1=a_1
a_1a_2a_4a_3z$. By Lemma~\ref{ij}, we have
\begin{eqnarray*}
a_1a_2a_4a_3za_1&=&a_1a_2za_4a_3a_1=a_1a_2z(a_1a_4a_3)\\
&=&a_1a_2a_1a_4za_3=a_1(a_1a_4a_2)za_3\\
&=&a_1a_1za_4a_2a_3=a_1a_1(a_2a_4a_3a_1)a_4a_2a_3\\
&=&a_1a_1a_2a_4a_3z.
\end{eqnarray*}
\end{proof}

\begin{lemma}\label{squares}
Let $z=a_1a_2\cdots a_n\in M$.
\begin{itemize}
\item[(i)] If $n\geq 6$ is even then
$a_i^2a_j (a_ka_la_rz)=a_ja_i^2 (a_ka_la_rz)$ and
 $a_ia_ja_ia_j (a_ka_la_rz)=a_ja_ia_ja_i
(a_ka_la_rz)$, for all $1\leq i,j\leq n$ and for any
three different
 integers $1\leq k,l,r\leq n$.
\item[(ii)] If $n\geq 4$, then $a_i^2a_jz^2=a_ja_i^2z^2$ and
$a_ia_ja_ia_jz^2=a_ja_ia_ja_iz^2$, for all $1\leq
i,j\leq n$.
\end{itemize}
\end{lemma}

\begin{proof}
$(i)$ Suppose that $n\geq 6$ is even. Applying
Lemma \ref{ij} several times, we get
\begin{eqnarray*}
a_{1}a_{1}a_{2} a_{1}a_{2}a_{3}z &=&
a_{1}a_{1}a_{2} (a_{3}a_{1}a_{2}z) =
a_{1}a_{1}a_{2} (z a_{3}a_{2}a_{1})\\ &=& a_{1}
(z a_{1}a_{2}) a_{3}a_{2}a_{1} = a_{1} (z
a_{2}a_{3} a_{1})a_{2}a_{1}\\ &=& a_{1} (a_{2}
a_{3} z) a_{1} a_{2}a_{1} = (z a_{2}a_{1}a_{3})
a_{1}a_{2} a_{1}\\ &=& (a_{2}a_{1}z)
a_{3}a_{1}a_{2}a_{1} = a_{2}a_{1} (z a_{1}
a_{2}a_{3})a_{1}\\ &=& a_{2}a_{1}(a_{1}a_{2}
a_{3}a_{1}z) = a_{2}a_{1}a_{1}(a_{1}a_{2}a_{3}z),
\end{eqnarray*}
\begin{eqnarray*}
a_{1}a_{2}a_{1}a_{2} a_{1}a_{2}a_{3}z &=&
a_{1}a_{2}a_{1}a_{2} (z a_{3}a_{2}a_{1}) =
a_{1}a_{2} (z a_{1}a_{2}) a_{3}a_{2}a_{1}\\ &=&
a_{1}a_{2} (a_{3}a_{2} a_{1}z)a_{2}a_{1} = a_{1}
a_{2} a_{3} (za_{2} a_{1}) a_{2}a_{1}\\ &=& (z
a_{2}a_{1}a_{3})a_{2} a_{1}a_{2} a_{1} =
(a_{2}a_{1}z) a_{3}a_{2}a_{1}a_{2}a_{1}\\ &=&
a_{2}a_{1} (z a_{2} a_{1}a_{3})a_{2}a_{1} =
a_{2}a_{1}(a_{2}a_{1}z) a_{3} a_{2}a_{1}\\ &=&
a_{2}a_{1}a_{2}a_{1}(a_{1}a_{2}a_{3}z)
\end{eqnarray*}
and, for every $i\in \{1,2,\dots ,n\}\setminus \{
3,4\}$,
\begin{eqnarray*}
a_{1}a_{1}a_{2} a_{i}a_{3}a_{4}z &=&
a_{1}a_{1}a_{2} (a_{3}a_{4}a_{i}z) =
a_{1}a_{1}a_{2} a_{3} (z a_{4}a_{i})\\ &=& a_{1}
(z a_{2}a_{1}a_{3})a_{4}a_{i} = a_{1} (z
a_{3}a_{2} a_{1})a_{4}a_{i}\\ &=& a_{1} (a_{3}
a_{2} z) a_{1} a_{4}a_{i} = (a_{2}a_{1}a_{3} z)
a_{1}a_{4} a_{i}\\ &=& a_{2}a_{1}a_{3}(a_{1}a_{4}
z) a_{i} = a_{2}a_{1} (a_{1} a_{4}a_{3} z)
a_{i}\\ &=& a_{2}a_{1}a_{1} (z a_{4} a_{3}) a_{i}
= a_{2}a_{1}a_{1}(a_{i}a_{3}a_{4}z),
\end{eqnarray*}
\begin{eqnarray*}
a_{1}a_{2}a_{1}a_{2} a_{i}a_{3}a_{4}z &=&
a_{1}a_{2}a_{1}a_{2} (za_{4}a_{3}a_{i}) =
a_{1}a_{2}(za_{1}a_{2})a_{4}a_{3}a_{i}\\ &=&
a_{1}a_{2} (a_{4}a_{2}a_{1}z)a_{3}a_{i} = a_{1}
a_{2}a_{4} (z a_{2}a_{1})a_{3}a_{i}\\ &=& (za_{2}
a_{1} a_{4}) a_{2}a_{1} a_{3}a_{i} =
(a_{2}a_{1}z)a_{4}a_{2} a_{1}a_{3} a_{i}\\ &=&
a_{2}a_{1}(z a_{2}a_{1}a_{4})a_{3} a_{i} =
a_{2}a_{1} (a_{2}a_{1}z) a_{4}a_{3}  a_{i}\\ &=&
a_{2}a_{1}a_{2}a_{1}(a_{i}a_{3}a_{4}z).
\end{eqnarray*}
Hence,  in each case applying an appropriate
$\sigma\in \Alt_n$ and using Lemma~\ref{ij}, we
obtain $a_i^2a_j (a_ka_la_rz)=a_ja_i^2
(a_ka_la_rz)$ and $a_ia_ja_ia_j
(a_ka_la_rz)=a_ja_ia_ja_i (a_ka_la_rz)$, for all
$1\leq i,j\leq n$ and for any three different
integers $1\leq k,l,r\leq n$.

$(ii)$  Suppose that $n$ is odd. Let
$z'=a_{6}a_{7}\cdots a_{n}$ (so $z'$ is the
identity element if $n=5$). Since $z$ is central
in $M$, by Lemma \ref{ij}, we have
\begin{eqnarray*}
a_{1}a_{1}a_{2}z^2&=&a_{1}a_{1}a_{2}(a_{3}a_{4}a_{5}a_{1}a_{2})z'z
=a_{1}(a_{2}a_{3}a_{1})a_{4}a_{5}a_{1}a_{2}z'z\\
&=&(a_{2}a_{3}a_{1})(a_{4}a_{5}a_{1})a_{1}a_{2}z'z
=a_{2}(a_{1}a_{4}a_{3})a_{5}a_{1}a_{1}a_{2}z'z\\
&=&a_{2}a_{1}a_{4}(a_{1}a_{3}a_{5})a_{1}a_{2}z'z
=a_{2}a_{1}(a_{1}a_{3}a_{4})a_{5}a_{1}a_{2}z'z\\
&=&a_{2}a_{1}a_{1}z^2
\end{eqnarray*}
and
\begin{eqnarray*}
a_1a_2a_1a_2z^2&=&
a_1a_2a_1a_2(a_3a_5a_4a_2a_1)z'z =
a_1a_2(a_3a_1a_2)a_5a_4a_2a_1z'z\\ &=&
a_1a_2a_3(a_5a_1a_2)a_4a_2a_1z'z =
a_1(a_5a_2a_3)(a_2a_4a_1)a_2a_1z'z \\ &=&
(a_2a_1a_5)a_3a_2a_4a_1a_2a_1z'z =
a_2a_1(a_2a_5a_3)(a_1a_2a_4)a_1z'z\\ &=&
a_2a_1a_2(a_1a_5a_3)a_2a_4a_1z'z=a_2a_1a_2a_1z^2.
\end{eqnarray*}

For $n=4$, we have
\begin{eqnarray*} a_1a_1a_2z^2&=&a_1z^2a_1a_2=a_1 (a_2a_3a_1a_4)( a_2a_3a_1a_4) a_1a_2\\
&=& a_1a_2a_3(a_4a_2a_1a_3) a_1a_4a_1a_2 =
(a_2a_1a_4a_3)a_2a_1a_3a_1a_4a_1a_2\\
&=&a_2a_1(a_1a_3a_4a_2)a_3a_1a_4a_1a_2 =
a_2a_1a_1a_3a_4 (a_1a_2a_3a_4) a_1a_2 \\ &=&
a_2a_1a_1z^2
\end{eqnarray*}
and, by Lemmas \ref{ij}, \ref{z2central} and
\ref{central4},
\begin{eqnarray*}
a_1a_2a_1a_2z^2&=& za_1a_2za_1a_2 =
a_2a_1(a_4a_3a_1a_2z)a_1a_2\\ &=&
a_2a_1(a_2a_1a_3a_4z)a_1a_2 =
a_2a_1a_2a_1a_3a_4a_1a_2z\\ &=& a_2a_1a_2a_1z^2.
\end{eqnarray*}

Hence, if $n$ is odd or $n=4$ and for $\sigma\in
\Alt_n$ we have that
$$a_{\sigma(1)}a_{\sigma(1)}a_{\sigma(2)}z^2=a_{\sigma(2)}a_{\sigma(1)}a_{\sigma(1)}z^2$$
and
$$a_{\sigma(1)}a_{\sigma(2)}a_{\sigma(1)}a_{\sigma(2)}z^2
=a_{\sigma(2)}a_{\sigma(1)}a_{\sigma(2)}a_{\sigma(1)}z^2.$$

So \begin{eqnarray*}
a_{i}^2a_{j}z^2=a_{j}a_{i}^2z^2\quad\mbox{and}\quad
a_{i}a_{j}a_{i}a_{j}z^2=a_{j}a_{i}a_{j}a_{i}z^2
\end{eqnarray*}
for all $1\leq i,j\leq n$, and $(ii)$ follows.

Suppose that $n\geq 6$ is even. By
Lemma~\ref{ij}, $z^2=a_2a_1a_3za_4a_5\cdots a_n$.
Consequently, by $(i)$, we get that
\begin{eqnarray*}
a_{i}^2a_{j}z^2=a_{j}a_{i}^2z^2\quad\mbox{and}\quad
a_{i}a_{j}a_{i}a_{j}z^2=a_{j}a_{i}a_{j}a_{i}z^2
\end{eqnarray*}
for all $1\leq i,j\leq n$, as desired.
\end{proof}

\begin{lemma}\label{normalf}
Let $z=a_1a_2\cdots a_n\in M$. For $1\leq i<j\leq
n$, let $F_{ij}=\langle a_i,a_j\rangle$. Then
\begin{itemize}
\item[(i)]
The elements in $z^2F_{ij}$ are of the form
\begin{eqnarray*}\label{Fij}
&&z^2a_i^{2n_1}a_j^{2n_2}w,
\end{eqnarray*}
where $w\in \{ 1, a_{i},a_{j}, a_{i}a_{j}, a_{j}a_{i},
a_{i}a_{j}a_{i}, a_{j}a_{i}a_{j}, a_{i}a_{j}a_{i}a_{j} \}$ and
$n_1,n_2$ are nonnegative integers.
\item[(ii)] The elements in $z^2(M\setminus
\bigcup_{1\leq i<j\leq n}F_{ij})$ are of the form
\begin{eqnarray}\label{Fijk}
z^2s_{i_1i_2}a_{i_2+1}^{m_1}a_{i_2+2}^{m_2}\cdots
a_n^{m_{n-i_2}},
\end{eqnarray}
where $i_1<i_2<n$, $s_{i_1i_2}\in
F_{i_1i_2}\setminus (\langle
a_{i_1}\rangle\cup\langle a_{i_2}\rangle)$,
$m_1,m_2,\dots, m_{n-i_2}$ are nonnegative
integers, and $\sum_{j=1}^{n-i_2}m_j>0$.
\end{itemize}
\end{lemma}

\begin{proof}
$(i)$ By Lemma~\ref{z2central}, $z^2$ is central
in $M$. Now, by Lemma~\ref{squares}, we have
\begin{eqnarray*}
a_ia_ja_ia_ja_iz^2&=& (a_ja_ia_ja_i)a_iz^2,\\
a_ia_ja_ia_ja_ia_jz^2&=& (a_ja_ia_ja_i)a_ia_jz^2=
a_ja_ia_ja_ja_i^2z^2,\\
a_ia_ja_ia_ja_ia_ja_iz^2&=&
(a_ja_ia_ja_i)a_ia_ja_iz^2 =  a_ja_i^4a_j^2z^2,\\
(a_ia_j)^4z^2&=&(a_ja_i)^2(a_ia_j)^2z^2 =
a_i^4a_j^4z^2.
\end{eqnarray*}
Therefore
\begin{eqnarray*}
(a_ia_j)^2a_iz^2&=&a_ja_ia_ja_i^2z^2, \\
(a_ia_j)^3z^2&=&a_ja_ia_j^2a_i^2z^2=a_ja_ia_i^2a_j^2z^2,\\
(a_ia_j)^3a_iz^2&=&a_ja_i^4a_j^2z^2=a_ja_j^2a_i^4z^2,\\
(a_ia_j)^4z^2&=&a_i^4a_j^4z^2=a_j^4a_i^4z^2,
\end{eqnarray*}
for all $1\leq i,j\leq n$. The above easily implies that the
elements in $z^2F_{ij}$ are of the form
\begin{eqnarray*}
&&z^2a_i^{2n_1}a_j^{2n_2}w,
\end{eqnarray*}
where $w\in \{ 1, a_{i},a_{j}, a_{i}a_{j}, a_{j}a_{i},
a_{i}a_{j}a_{i}, a_{j}a_{i}a_{j}, a_{i}a_{j}a_{i}a_{j} \}$ and
$n_1,n_2$ are nonnegative integers.

$(ii)$ Let $1\leq i<j\leq n$ and $k\in \{ 1,\dots
,n\}\setminus \{ i,j\}$. Then, by~(\ref{Fij}) and
Lemmas~\ref{ij}, \ref{z2central}
and~\ref{squares}, it is easy to see that
\begin{eqnarray}\label{kFij}
a_kz^2(F_{ij}\setminus(\langle
a_i\rangle\cup\langle a_j\rangle))=
z^2(F_{ij}\setminus(\langle a_i\rangle\cup\langle
a_j\rangle))a_k,
\end{eqnarray}
for all $n\geq 4$.

Let $s\in M\setminus \bigcup_{1\leq i<j\leq
n}F_{ij}$. Then $s=a_{j_1}a_{j_2}\cdots a_{j_k}$,
where $\{ j_1,\dots ,j_k\}$ is a subset of $\{
1,\dots ,n\}$ of cardinality $\geq 3$. We shall
prove that $z^2s$ is of the form~(\ref{Fijk}) by
induction on the total degree $k\geq 3$ of $s$.
For $k=3$, we have that $j_1,j_2,j_3$ are three
different elements and, by Lemmas~\ref{ij}
and~\ref{z2central},
$$z^2a_{j_1}a_{j_2}a_{j_3}=z^2a_{j_2}a_{j_3}a_{j_1}=z^2a_{j_3}a_{j_1}a_{j_2},$$
thus the result follows in this case.

Suppose that $k>3$ and that the result is true
for all elements in $M\setminus \bigcup_{1\leq
i<j\leq n}F_{ij}$ of total degree less than $k$.
Then either $a_{j_2}\cdots a_{j_k}\in
F_{i_1i_2}\setminus(\langle
a_{i_1}\rangle\cup\langle a_{i_2}\rangle)$, for
some $i_1<i_2$, or $a_{j_2}\cdots a_{j_k}\in
M\setminus \bigcup_{1\leq i<j\leq n}F_{ij}$. By
the induction hypothesis $$z^2a_{j_2}\cdots
a_{j_{k}}=z^2s_{i_1i_2}a_{i_2+1}^{m_1}a_{i_2+2}^{m_2}\cdots
a_n^{m_{n-i_2}},$$ where $i_1<i_2\leq n$,
$s_{i_1i_2}\in F_{i_1i_2}\setminus(\langle
a_{i_1}\rangle\cup\langle a_{i_2}\rangle)$ and
$m_1,\dots ,m_{n-i_2}\geq 0$. We may assume that
$j_1\notin\{ i_1,i_2\}$ and, by~(\ref{kFij}),
$$z^2s=z^2s_{i_1i_2}a_{j_1}a_{i_2+1}^{m_1}a_{i_2+2}^{m_2}\cdots
a_n^{m_{n-i_2}}.$$
 Suppose that $j_1<i_2$. In this case, by Lemmas~\ref{ij} and~\ref{z2central}, we get
 \begin{eqnarray*}
a_{i_2}a_{i_1}a_{i_2}a_{j_1}z^2&=&a_{i_2}(a_{i_2}a_{j_1}a_{i_1})z^2
= a_{i_2}^2a_{j_1}a_{i_1}z^2,\\
a_{i_1}a_{i_2}a_{i_1}a_{j_1}z^2&=&a_{i_1}(a_{i_1}a_{j_1}a_{i_2})z^2
=a_{i_1}^2a_{j_1}a_{i_2}z^2,\\
(a_{i_1}a_{i_2})^2a_{j_1}z^2&=&a_{i_1}a_{i_2}(a_{i_2}a_{j_1}a_{i_1})z^2
=(a_{i_2}^2)a_{i_1}a_{j_1}a_{i_1}z^2, \\
(a_{i_1}a_{i_2})^2a_{i_1}a_{j_1}z^2&=&a_{i_1}a_{i_2}(a_{i_1}^2a_{j_1}a_{i_2})z^2
= (a_{i_1}^2)a_{i_1}a_{i_2}a_{j_1}a_{i_2}z^2\\
 &=& a_{i_1}^2(a_{j_1}a_{i_1}a_{i_2})a_{i_2}z^2=a_{i_1}^2(a_{i_2}^2)a_{j_1}a_{i_1}z^2,\\
(a_{i_1}a_{i_2})^3a_{j_1}z^2&=&a_{i_1}a_{i_2}(a_{i_2}^2a_{i_1}a_{j_1}a_{i_1})z^2
=(a_{i_2}^2)a_{i_1}a_{i_2}a_{i_1}a_{j_1}a_{i_1}z^2\\
&=&
a_{i_2}^2a_{i_1}(a_{i_1}a_{j_1}a_{i_2})a_{i_1}z^2
= (a_{i_1}^2)a_{i_2}^2a_{j_1}a_{i_2}a_{i_1}z^2\\
&=&
a_{i_1}^2a_{i_2}^2(a_{i_1}a_{j_1}a_{i_2})z^2,\\
(a_{i_1}a_{i_2})^3a_{i_1}a_{j_1}z^2&=&a_{i_1}a_{i_2}(a_{i_1}^2a_{i_2}^2a_{j_1}a_{i_1})z^2
=
(a_{i_1}^2a_{i_2}^2)a_{i_1}a_{i_2}a_{j_1}a_{i_1}z^2\\
&=&
a_{i_1}^2a_{i_2}^2a_{i_1}(a_{j_1}a_{i_1}a_{i_2})z^2.
\end{eqnarray*}
By Lemma~\ref{squares} and~(\ref{Fij}), the
result follows, in this case.

Suppose that $j_1>i_2$. By Lemma~\ref{squares},
$s_{i_1i_2}\in (F_{i_1i_2}a_{i_1}a_{i_2})\cup
(F_{i_1i_2}a_{i_2}a_{i_1})$. Note that if
$i_2<k<j_1$, then
$$za_{i_1}a_{i_2}a_{j_1}a_k=za_{i_2}a_{i_1}a_ka_{j_1}
\mbox{ and }
za_{i_2}a_{i_1}a_{j_1}a_k=za_{i_1}a_{i_2}a_ka_{j_1}.$$ Therefore
(\ref{kFij}) implies that $z^2s$ is of the form~(\ref{Fijk}). Thus
the result follows by induction.
\end{proof}

\begin{lemma}\label{normalf2}
Suppose that $n\geq 6$ is even. Let
$z=a_1a_2\cdots a_n\in M$. For $1\leq i<j\leq n$,
let $F_{ij}=\langle a_i,a_j\rangle$. Let $k,l,r$
be three different integers such that $1\leq
k,l,r\leq n$. Then
\begin{itemize}
\item[(i)] $(a_ka_la_rz)a_i=a_i(a_ka_la_rz)$, for all
 $i\in\{ 1,2,\dots ,n\}\setminus\{ k,l,r\}$.
\item[(ii)] $(a_ka_la_rz)a_i=a_i(a_la_ka_rz)$, for all $i\in\{ k,l,r\}$.
\item[(iii)]
The elements in $a_ka_la_rzF_{ij}$ are of the form
\begin{eqnarray*}
&&a_ka_la_rza_i^{2n_1}a_j^{2n_2}w,
\end{eqnarray*}
where $w\in \{ 1, a_{i},a_{j}, a_{i}a_{j}, a_{j}a_{i},
a_{i}a_{j}a_{i}, a_{j}a_{i}a_{j}, a_{i}a_{j}a_{i}a_{j} \}$ and
$n_1,n_2$ are nonnegative integers.
\item[(iv)] The elements in $a_ka_la_rz(M\setminus
\bigcup_{1\leq i<j\leq n}F_{ij})$ are of the form
\begin{eqnarray*}
a_ka_la_rzs_{i_1i_2}a_{i_2+1}^{m_1}a_{i_2+2}^{m_2}\cdots
a_n^{m_{n-i_2}},
\end{eqnarray*}
where $i_1<i_2<n$, $s_{i_1i_2}\in
F_{i_1i_2}\setminus (\langle
a_{i_1}\rangle\cup\langle a_{i_2}\rangle)$,
$m_1,m_2,\dots, m_{n-i_2}$ are nonnegative
integers, and $\sum_{j=1}^{n-i_2}m_j>0$.
\end{itemize}
\end{lemma}

\begin{proof}
$(i)$ Let $i\in\{ 1,2,\dots ,n\}\setminus\{
k,l,r\}$. By Lemma~\ref{ij}, we have
\begin{eqnarray*}
(a_ka_la_rz)a_i &=&
a_k(za_la_r)a_i=a_kz(a_ia_la_r)=a_k(a_ia_lz)a_r\\
&=&
(a_ia_la_k)za_r=a_i(za_la_k)a_r=a_i(a_ka_la_rz).
\end{eqnarray*}

$(ii)$ Let $i\in\{ k,l,r\}$. By Lemma~\ref{ij},
we may assume that $i=k$, and we have
\begin{eqnarray*}
(a_ka_la_rz)a_k =
a_k(za_la_r)a_k=a_kz(a_ka_la_r)=a_k(a_la_ka_rz).
\end{eqnarray*}

The proof of $(iii)$ and $(iv)$ is similar to the
proof of Lemma~\ref{normalf}.  Namely, it is
obtained by using $(i)$ and $(ii)$ in place of
the fact that $z^2$ is central and using
Lemma~\ref{squares}$(i)$ in place of
Lemma~\ref{squares}$(ii)$.
\end{proof}

\begin{lemma}\label{even1}
Suppose that $n\geq 6$ is even. Then
$$\bigcup_{1\leq r\leq n}(Mz\cap Mza_r) =
\bigcup_{1\leq i<j<k\leq n}(Ma_ia_ja_kz\cup
Ma_ja_ia_kz).$$
\end{lemma}

\begin{proof}
By Lemma \ref{ij}, we have that
\begin{eqnarray*}
a_1a_2a_3z= za_2 a_1a_3=a_2a_1za_3.
\end{eqnarray*}
Note that if $1\leq i,j,k\leq n$ are three
different integers then, since $n\geq 6$, there
exists $\sigma\in \Alt_n$ such that
$\sigma(1)=i$, $\sigma(2)=j$ and $\sigma(3)=k$.
Therefore
\begin{eqnarray}\label{MzMzr}
a_ia_ja_kz\in \bigcup_{1\leq r\leq n}(Mz\cap
Mza_r),
\end{eqnarray}
for all different $1\leq i,j,k\leq n$.

Suppose that $\bigcup_{1\leq r\leq n}(Mz\cap
Mza_r)\not\subseteq \bigcup_{1\leq i<j<k\leq
n}(Ma_ia_ja_kz\cup Ma_ja_ia_kz)$. Let $s\in
\bigcup_{1\leq r\leq n}(Mz\cap Mza_r)\setminus
\bigcup_{1\leq i<j<k\leq n}(Ma_ia_ja_kz\cup
Ma_ja_ia_kz)$ be an element of minimal length.
There exist $1\leq r\leq n$, $s'=a_{j_1}\cdots
a_{j_{k-1}}\in M$ and $s''=a_{i_1}a_{i_2}\cdots
a_{i_k}$ such that $s=s'za_r=s''z$. Thus there
exist $w_1,w_2,\dots ,w_m$ in the free monoid
$\FM_{n}$ on $\{ a_1,\dots ,a_n\}$, such that
$w_1=a_{j_1}\cdots a_{j_{k-1}}a_1a_2\cdots
a_na_r$, $w_m=a_{i_1}\cdots a_{i_{k}}a_1a_2\cdots
a_n$ and
$w_i=w_{1,i}w_{2,i}w_{3,i}=w'_{1,i}w'_{2,i}w'_{3,i}$,
where $w_{2,i}$ and $w'_{2,i}$ represent the
element $z$ in $M$ for all $i=1,\dots ,m$, and
$w_{1,j}=w'_{1,j+1}$ and $w_{3,j}=w'_{3,j+1}$,
for all $j=1,\dots ,m-1$.

Let $g\colon\{ 1,2,\dots ,m\}\times \{ 1,2,\dots
,n+k\}\longrightarrow \{ 1,2,\dots ,n\}$ be such
that $w_i=a_{g(i,1)}a_{g(i,2)}\cdots
a_{g(i,n+k)}$ for all $i=1,\dots ,m$. Let $t$ be
the least positive integer such that
$a_{g(t,k+1)}a_{g(t,k+2)}\cdots a_{g(t,n+k)}$
represents $z$ in $M$. Since $n$ is even, $t>1$
and $g(i,n+k)=g(t,n+k)=r$, for all $i=1,\dots
,t$. Hence $$a_{g(1,1)}a_{g(1,2)}\cdots
a_{g(1,n+k-1)},\;\dots ,\;
a_{g(t-1,1)}a_{g(t-1,2)}\cdots a_{g(t-1,n+k-1)}$$
represent the same element in $M$. Furthermore,
the length of $w_{3,t-1}$ is less than $n$ and
greater than $0$.

Suppose that $w_{3,t-1}=a_{r}$. In this case,
$w'_{2,t}a_{r}=a_{g(t,k)}\cdots a_{g(t,n+k)}$ and
$w_{2,t-1}a_{r}$ represent the same element in
$M$, but in $M$ we have that $za_r\neq
a_{g(t,k)}z$, a contradiction. Therefore the
length of $w_{3,t-1}$ is greater than $1$.

Suppose that $w_{3,t-1}=a_{g(t-1,n+k-1)}a_{r}$.
In this case, $w_{2,t-1}a_{g(t-1,n+k-1)}a_{r}$
and
$w'_{2,t}a_{g(t-1,n+k-1)}a_{r}=a_{g(t,k-1)}a_{g(t,k)}\cdots
a_{g(t,n+k)}$ represent the same element in $M$.
Since $$a_{g(1,1)}a_{g(1,2)}\cdots
a_{g(1,n+k-1)},\;\dots ,\;
a_{g(t-1,1)}a_{g(t-1,2)}\cdots a_{g(t-1,n+k-1)}$$
represent the same element in $M$, we have in $M$
that
\begin{eqnarray*}
a_{g(1,1)}\cdots
a_{g(1,k-1)}z&=&a_{g(1,1)}a_{g(1,2)}\cdots
a_{g(1,n+k-1)}\\
&=&a_{g(t-1,1)}a_{g(t-1,2)}\cdots
a_{g(t-1,n+k-1)}\\
&=&a_{g(t-1,1)}a_{g(t-1,2)}\cdots
a_{g(t-1,k-2)}za_{g(t-1,n+k-1)}.
\end{eqnarray*}
Thus $a_{g(1,1)}\cdots a_{g(1,k-1)}z\in Mz\cap
Mza_{g(t-1,n+k-1)}$. By the choice of $s$,  we
have that $$a_{g(1,1)}\cdots a_{g(1,k-1)}z\in
\bigcup_{1\leq i<j<k\leq n}(Ma_ia_ja_kz\cup
Ma_ja_ia_kz).$$ Since $s=a_{g(1,1)}\cdots
a_{g(1,k-1)}za_r$, by Lemma \ref{normalf2} $(i)$ and $(ii)$,
$$s\in
\bigcup_{1\leq i<j<k\leq n}(Ma_ia_ja_kz\cup
Ma_ja_ia_kz),$$ a contradiction. Therefore the
length of $w_{3,t-1}$ is greater than $2$.

Thus $w_{3,t-1}=a_{g(t-1,n+k-l)}\cdots
a_{g(t-1,n+k-1)}a_{r}$ for some $1<l<n$.  In this
case, $$w'_{2,t}a_{g(t-1,n+k-l)}\cdots
a_{g(t-1,n+k-1)}a_{r} =a_{g(t,k-l)}\cdots
a_{g(t,k-1)}a_{g(t,k)}\cdots a_{g(t,n+k)}.$$
Since $a_{g(t,k+1)}\cdots a_{g(t,n+k)}$
represents $z$ in $M$ and $l<n$, we have that
$g(t-1,n+k-l),\ldots , g(t-1,n+k-1), r$ are $l+1$
different integers. Hence $s\in
Mza_{g(t-1,n+k-l)}\cdots a_{g(t-1,n+k-1)}a_{r}$.
By Lemma~\ref{ij}, $$s\in \bigcup_{1\leq
i<j<k\leq n}(Ma_ia_ja_kz\cup Ma_ja_ia_kz),$$ a
contradiction. Therefore $$\bigcup_{1\leq r\leq
n}(Mz\cap Mza_r)\subseteq \bigcup_{1\leq
i<j<k\leq n}(Ma_ia_ja_kz\cup Ma_ja_ia_kz).$$ By
(\ref{MzMzr}), the result follows.
\end{proof}

\begin{lemma}\label{preliminar}
Suppose that $n\geq 6$ is even. Let
$s=a_{j_1}a_{j_2}\cdots a_{j_m}\in M\setminus
MzM$ such that $$sz\notin\bigcup_{1\leq i<j<k\leq
n}(Ma_{i}a_{j}a_{k}z\cup Ma_{j}a_{i}a_{k}z).$$
Then, for $s_1,s_2\in M$, $sz=s_1zs_2$ implies
that $s_1s_2=s$.
\end{lemma}

\begin{proof}
Let $s_1,s_2\in M$ such that $sz=s_1zs_2$. Then,
by an easy degree argument, $s_1=a_{i_1}\cdots
a_{i_k}$ and $s_2=a_{i_{k+1}}\cdots a_{i_{m}}$
for some $k$ and some $a_{i_{1}},\ldots
,a_{i_{m}}$. Thus there exist $w_1,w_2,\dots
,w_t$ in the free monoid $\FM_{n}$ on $\{
a_1,\dots ,a_n\}$, such that
$w_i=w_{1,i}w_{2,i}w_{3,i}=w'_{1,i}w'_{2,i}w'_{3,i}$,
where $w_{2,i}$ and $w'_{2,i}$ represent the
element $z$ in $M$ for all $i=1,\dots ,t$,
$w_{1,j}=w'_{1,j+1}$ and $w_{3,j}=w'_{3,j+1}$,
for all $j=1,\dots ,t-1$, and
$w'_{1,1}=a_{j_1}a_{j_2}\cdots a_{j_m}$,
$w'_{3,1}=1$, $w_{1,t}=a_{i_1}\cdots a_{i_k}$ and
$w_{3,t}=a_{i_{k+1}}\cdots a_{i_{m}}$. Thus,
$w_{1}=a_{j_{1}}\cdots a_{j_{m}} w_{2,1}'$ and
$w_{t}=a_{i_{1}}\cdots a_{i_{k}} w_{2,t}
a_{i_{k+1}}\cdots a_{i_{m}}$. It is enough to
prove that $w_{1,i}w_{3,i}=a_{j_1}\cdots
a_{j_m}$, for all $i=1,\dots ,t$, by induction on
$t$.

If the two subwords $w_{2,1}$ and $w_{2,1}'$ of the word $w_1=
a_{j_1}a_{j_2}\cdots a_{j_m}w'_{2,1} = w_{1,1}w_{2,1}w_{3,1}$ do
not overlap, then $w_{2,1}$ is a subword of $a_{j_{1}}\cdots
a_{j_{m}}$, which is not possible because the latter represents
$s$ in $M$ and $s\notin MzM$. Therefore they overlap and hence the
degree of $w_{3,1}$ is less than $n$ and $w_{3,1}$ is a product of
distinct letters. Since $w_{1}$ represents $sz$ in $M$ and
$sz\notin\bigcup_{1\leq i<j<k\leq n}(Ma_{i}a_{j}a_{k}z\cup
Ma_{j}a_{i}a_{k}z)$, it follows that $w_{3,1}$ cannot have degree
$1$ by Lemma~\ref{even1} and it cannot have degree greater than
$2$ by Lemma~\ref{ij}. Hence, the degree of $w_{3,1}$ is $0$ or
$2$. In the former case, clearly
$w_{1,1}w_{3,1}=w_{1,1}=w'_{1,1}=a_{j_1}a_{j_2}\cdots a_{j_m}$.
Suppose that $w_{3,1}$ has degree $2$. From the equality of words
$a_{j_1}a_{j_2}\cdots a_{j_m}w'_{2,1}=w_{1,1}w_{2,1}w_{3,1}$ it
follows that $a_{j_{m-1}}a_{j_{m}}w_{2,1}'=w_{2,1}w_{3,1}$. Since
in $M$ $a_ia_jz=za_ka_l$ if and only if $i\neq j$ and
$(i,j)=(k,l)$, this implies that that
$w_{3,1}=a_{j_{m-1}}a_{j_{m}}$ and clearly
$w_{1,1}=a_{j_1}a_{j_2}\cdots a_{j_{m-2}}$. Hence, in both cases,
$w_{1,1}w_{3,1}=a_{j_1}\cdots a_{j_m}$. Suppose that $t>1$ and
$w_{1,i}w_{3,i}=a_{j_1}\cdots a_{j_m}$, for all $i=1,\dots ,t-1$.

We have that $w'_{1,t}=w_{1,t-1}=a_{j_1}\cdots
a_{j_q}$ and
$w'_{3,t}=w_{3,t-1}=a_{j_{q+1}}\cdots a_{j_m}$,
for some $0\leq q\leq m$. Hence
$w_t=a_{j_1}\cdots
a_{j_q}w'_{2,t}a_{j_{q+1}}\cdots
a_{j_m}=w_{1,t}w_{2,t}w_{3,t}$. Since
$a_{j_{1}}\cdots a_{j_{m}}\not \in MzM$ by the
hypothesis, as above we get that the subwords
$w_{2,t}$ and $w_{2,t}'$ of the word $w_{t}$ have
to overlap.

Let $r$ be the absolute value of the difference
of the lengths of the words $w_{1,t}$ and
$w_{1,t}'$. Then $r<n$. The equality of words
$w_{1,t}'w_{2,t}'w_{3,t}'=w_{1,t}w_{2,t}w_{3,t}$
implies that either $w_{2,t}'u'=uw_{2,t}$ or
$u'w_{2,t}'=w_{2,t}u$ for some words $u,u'$ of
length $r$. Then all the generators involved in
$u$ (and also in $u'$) are different. Since $w_t$
represents $sz$ in $M$ and
$sz\notin\bigcup_{1\leq i<j<k\leq
n}(Ma_{i}a_{j}a_{k}z\cup Ma_{j}a_{i}a_{k}z)$, by
Lemma~\ref{ij} we get that $r\leq 2$. If $r=1$
then we get $za_{p}=a_{p}z$ in $M$ for some $p$,
which is impossible since $n$ is even. Hence,
$r=0$ or $r=2$. As above we can see in the both
cases that $w_{1,t}w_{3,t}=a_{j_1}\cdots
a_{j_m}$. The result follows.
\end{proof}

Let $I=\{ s\in M\mid sM\subseteq Mz\}$ and $I'=\{
s\in M\mid Ms\subseteq zM\}$. Clearly, $I$ and
$I'$ are ideals of $M$. Let $I_1=\{ s\in Mz\mid
sa_i\in Mz,$ for all $i=1,2,\dots ,n \}$,
$I'_1=\{ s\in zM\mid a_is\in zM,$ for all
$i=1,2,\dots ,n \}$ and $$T=\bigcup_{1\leq
i<j<k\leq n}(Ma_ia_ja_kz\cup Ma_ja_ia_kz).$$

Let $u,u'$ be words in the free monoid $\FM_n$ on $\{
a_1,a_2,\ldots, ,a_n\}$. We say that $u'$ is a one step rewriting
of $u$ if there exist $u_1,u_2,u_3, u_{2}'$ in $\FM_n$ such that
$u=u_{1}u_2u_{3}, u'=u_{1}u_{2}'u_{3}$ and both $u_2, u_{2}'$
represent $z$ in $M_{n}$.

\begin{lemma}\label{even2}
Suppose that $n\geq 6$ is even. Then
$I=I'=I_1=I'_1=T$.
\end{lemma}

\begin{proof}
By Lemma~\ref{normalf2} $(i)$ and $(ii)$, we have
that $T$ is an ideal of $M$. Hence $T\subseteq
I$. Suppose that these two ideals are different.
Let $s\in I\setminus T$. Since $s\in I$, there
exists $s'\in M$ such that $s=s'z$. We consider
two cases.

{\em Case $1$:} $s'\in MzM$.

Then let $s''\in M$ be an element of minimal degree such that
$s'\in Mzs''$. Thus there exists $t\in M$ such that $s'=tzs''$.
Since $s=tzs''z\notin T$, we have that $s''$ has degree greater
than or equal to $2$. By Lemma~\ref{ij} and the choice of $s''$,
there exists $i$ such that $s''\in a_i^2M\setminus MzM$, so
$s''=a_i^2a_{j_1}\cdots a_{j_m}$ for some $m\geq 0$. Let $s_1\in
M$ be an element of minimal degree such that $s''z\in s_1zM$. Then
there exists $s_2\in M$ such that $s''z=s_1zs_2$. Since
$s=tzs''z=tzs_1zs_2\notin T$, it follows that $s_1$ has degree
greater than or equal to $2$. Since $s''\notin MzM$ and
$s''z=s_1zs_2\notin T$, by Lemma~\ref{preliminar}, we have that
$s''=s_1s_2$. Let $$h=\left\{\begin{array}{ll} i&\quad\mbox{if }
m=0\\ j_m&\quad\mbox{if } m>0 .\end{array}\right.$$ Note that
$sa_h=tzs''za_h=tzs_1zs_2a_h$ and $s'=tzs''=tzs_1s_2$. By the
choice of $s''$, if we rewrite $tzs_1$, this has to be of the form
$t's_1$ for some $t'$. Note also that all rewritings of $s_1zs_2$
must be of the form: $s'_1zs'_2$ where $s'_1s'_2=s''$. Now we look
at the one step rewritings of $s'_1zs'_2a_h$. If $s'_2=1$ then
every one step rewriting of $s'_1zs'_2a_h$ cannot affect the last
$a_h$ because $n$ is even. If $s'_2\neq 1$ then every one step
rewriting of $s'_1zs'_2a_h$ cannot affect the last $a_h$ because
the last generator in $s_2'$ is also $a_h$. By the choice of
$s_{1}$ all possible rewritings of $s_{1}zs_{2}a_{h}$ are of the
form $s_{1}wa_{h}$ for some $w$. But this shows that
$sa_h=tzs''za_h=tzs_1zs_2a_h\notin Mz$, a contradiction since
$s\in I$.

{\em Case $2$:} $s'\notin MzM$.

Then $s'=a_{j_{1}}\cdots a_{j_{m}}$ for some $m\geq 0$. If $m=0$
then $s=s'z=z$ and $sa_{1}=za_{1}\notin Mz$, a contradiction.
Hence $m>0$. Suppose $s_{1},s_{2}\in M$ are such that
$s=s'z=s_{1}zs_{2}$. By Lemma~\ref{preliminar} we get
$s'=s_{1}s_{2}$. Using a similar argument as in Case~1, one can
show that all the rewritings of $sa_{j_{m}}$ cannot affect the
last generator $a_{j_{m}}$. Therefore
$sa_{j_{m}}=s'za_{j_{m}}\notin Mz$, a contradiction since $s\in
I$. Therefore $I=T$.

Clearly, we have $I\subseteq I_1$. Let $s\in I_1$
and let $t\in M\setminus \{ 1\}$. Then
$t=a_{r}t'$ for some $1\leq r \leq n$ and $t'\in
M$. Since $sa_{r}\in Mz\cap Mza_{r}$, by
Lemma~\ref{even1} it follows that $st=sa_{r}t'\in
Tt'\subseteq T\subseteq Mz$. Therefore $s\in I$
and so $I=I_1$.

By Lemma~\ref{ij} and Lemma~\ref{normalf2} $(i)$
and $(ii)$, $$\bigcup_{1\leq i<j<k\leq
n}(za_ia_ja_kM\cup za_ja_ia_kM)=\bigcup_{1\leq
i<j<k\leq n}(Ma_ia_ja_kz\cup Ma_ja_ia_kz)=T.$$
Thus, by symmetry, $$I=I'=I_1=I'_1=T.$$
\end{proof}

\section{Proof of Theorem~\ref{main-alt}} \label{altern}

In this section we prove our main result,
Theorem~\ref{main-alt}. So again,  $n\geq 4$,
$M=S_n(\Alt_n)$ and $G=G_n(\Alt_n)$.

Recall that $\rho'$ is the binary relation on $M$, defined by
$s\rho' t$ if and only if there exists a nonnegative integer $i$
such that $sz^{i}=tz^{i}$. By Lemma~\ref{z2central}, $z^2$ is
central in $M$. By Lemma~\ref{canc}, $\rho'=\rho$ is the least
cancellative congruence on $M$.

{\it Proof of $(i)$.}  Let $\{ i,j,k\}$ be a
subset of $\{ 1,2,\dots ,n\}$ of cardinality
three. By Lemma~\ref{ij},  in $G$, we have
\begin{eqnarray}\label{ijk}
a_{i}a_{j}a_{k}=a_{j}a_{k}a_{i}=a_{k}a_{i}a_{j}.
\end{eqnarray}
By Lemma~\ref{squares}, in $G$, we also have
$a_ia_ja_ia_j=a_ja_ia_ja_i$ and
$a_{i}^{2}a_{j}=a_{j}a_{i}^{2}$, for all $1\leq
i,j\leq n$. Therefore
$(a_ia_ja_i^{-1}a_j^{-1})^2=1$.

Let $\tau\in \Sym_n\setminus \Alt_n$. By
(\ref{ijk}), it is easy to see that in $G$ we
have that
\begin{eqnarray*}
(2\; 3)(a_1a_2\cdots a_n)= a_1a_3a_2a_4\dots
a_n=a_{\tau(1)}a_{\tau(2)}\cdots a_{\tau(n)}.
\end{eqnarray*}
Hence we have the following presentations of the
group $G$.
\begin{eqnarray*}G&=& \gr( a_1,\dots ,a_n\mid
a_1a_2\cdots a_n=a_{\sigma(1)}a_{\sigma(2)}\cdots
a_{\sigma(n)}, \sigma\in \Alt_n)\\ &=& \gr(
a_1,\dots ,a_n\mid a_1a_2\cdots
a_n=a_{\sigma(1)}a_{\sigma(2)}\cdots
a_{\sigma(n)},\\ &&\hphantom{gra}
a_1a_3a_2a_4\cdots
a_n=a_{\tau(1)}a_{\tau(2)}\cdots a_{\tau(n)},\,
\sigma\in \Alt_n,\,  \tau\in \Sym_n\setminus
\Alt_n ).
\end{eqnarray*}
Note that, by (\ref{ijk}),
$a_i(a_1a_2a_1^{-1}a_2^{-1})a_i^{-1}=a_1a_2a_1^{-1}a_2^{-1}$,
for all $2<i\leq n$. Furthermore
\begin{eqnarray*}
a_1(a_1a_2a_1^{-1}a_2^{-1})a_1^{-1}&=&a_1(a_1a_2a_3)(a_3^{-1}a_1^{-1}a_2^{-1})a_1^{-1}\\
&=&a_1(a_2a_3a_1)(a_1^{-1}a_2^{-1}a_3^{-1})a_1^{-1}\quad\mbox{by
(\ref{ijk})}\\
&=&a_1a_2a_3a_2^{-1}a_3^{-1}a_1^{-1}\\
&=&a_1a_2a_3(a_3^{-1}a_1^{-1}a_2^{-1})\quad\mbox{by
(\ref{ijk})}\\ &=&a_1a_2a_1^{-1}a_2^{-1}
\end{eqnarray*}
and
\begin{eqnarray*}
a_2(a_1a_2a_1^{-1}a_2^{-1})a_2^{-1}&=&a_2(a_1a_2a_3)(a_3^{-1}a_1^{-1}a_2^{-1})a_2^{-1}\\
&=&a_2(a_3a_1a_2)(a_2^{-1}a_3^{-1}a_1^{-1})a_2^{-1}\quad\mbox{by
(\ref{ijk})}\\
&=&a_2a_3a_1a_3^{-1}a_1^{-1}a_2^{-1}\\
&=&(a_1a_2a_3)a_3^{-1}a_1^{-1}a_2^{-1}\quad\mbox{by
(\ref{ijk})}\\ &=&a_1a_2a_1^{-1}a_2^{-1} .
\end{eqnarray*}
Therefore $a_1a_2a_1^{-1}a_2^{-1}$ is a central
element of order at most $2$ in $G$. Let $C$ be
the central subgroup $C=\{
1,a_1a_2a_1^{-1}a_2^{-1}\}$. Then $G/C$ has the
following presentations.
\begin{eqnarray*}G/C&=& \gr( b_1,\dots ,b_n\mid b_1b_2=b_2b_1,\;
b_1b_2\cdots b_n=b_{\sigma(1)}b_{\sigma(2)}\cdots
b_{\sigma(n)},\\ &&\hphantom{\gr(}
b_1b_3b_2b_4\cdots
b_n=b_{\tau(1)}b_{\tau(2)}\cdots b_{\tau(n)},\,
\sigma\in \Alt_n,\, \tau\in \Sym_n\setminus
\Alt_n )\\ &=& \gr( b_1,\dots ,b_n\mid
b_1b_2\cdots b_n=b_{\sigma(1)}b_{\sigma(2)}\cdots
b_{\sigma(n)},\,  \sigma\in \Sym_n ).
\end{eqnarray*}
Hence $G/C$ is a free abelian group of rank $n$
and, since $C=G'$, in $G$ we have
\begin{eqnarray*}
a_1a_2a_1^{-1}a_2^{-1}=a_ia_ja_i^{-1}a_j^{-1},
\end{eqnarray*}
for all $i\neq j$, because there exists
$\sigma\in \Alt_n$ such that $\sigma(1)=i$ and
$\sigma(2)=j$.

We now show that $a_1a_2a_1^{-1}a_2^{-1}\neq 1$.
For this, let $\FM_{n}$ be the free monoid on the
set $\{ a_{1},\dots,a_{n}\}$. We define the map
$f\colon \FM_{n}\longrightarrow \{ -1,1\}$ by
\begin{eqnarray}\label{f}
f(a_{i_{1}}\cdots
a_{i_{m}})=\prod_{\begin{array}{c}
\scriptstyle{1\leq j<k\leq m}\\
\scriptstyle{i_{j}\neq
i_{k}}\end{array}}\frac{i_{k}-i_{j}}{|i_{k}-i_{j}|}.
\end{eqnarray}
Note that if two words $w,w'\in \FM_{n}$
represent the same element in $M_n$ then
$f(w)=f(w')$. In particular, $$a_{1}a_{2}z^m\neq
a_{2}a_{1}z^m$$ in $M$, for all $m$. Now, by
Lemmas~\ref{z2central} and~\ref{canc}, we have
that $a_1a_2\neq a_2a_1$ in $G$, as desired.

As seen above, every $a_i^2$ is a central element
in $G$. Let $D$ be the central subgroup of $G$
generated by $a_1^2,\dots , a_n^2$. Now $G/D$ is
finite because $C$ is finite and $G/(CD)\cong
(\mathbb{Z}/2\mathbb{Z})^n$. Hence $(i)$ follows.

{\it Proof of $(ii)$.}  By $(i)$,  $K[G]$ is a
noetherian PI-algebra for any field $K$.
Furthermore, by \cite[Theorem~7.3.1]{passman}
$\mathcal{J}(K[G])\subseteq
\mathcal{J}(K[C])K[G]$. Thus, if K is a field of
characteristic $\neq 2$, then
$\mathcal{J}(K[G])=0$. If $K$ is a field of
characteristic $2$, then
$\mathcal{J}(K[G])=(1-a_1a_2a_1^{-1}a_2^{-1})K[G]$,
and $\mathcal{J}(K[G])^2=0$.

{\it Proof of $(iii)$.}  By
Lemma~\ref{z2central}, $z^2$ is central in $M$.
Note that two elements of different
forms~(\ref{Fij}) in $M$ are different in $G$ and
they thus also differ in $M$. Also, two elements
of $M$ that have different forms~(\ref{Fijk}) are
different in $G$ and thus also in $M$. Moreover,
elements in (\ref{Fij}) and (\ref{Fijk}) are
different. Since $z^{2}$ is central in $M$, it
follows from Lemma~\ref{canc} that every right
ideal of $M$ contains $z^{2k}$ for some positive
integer $k$. Therefore, if $sx=tx$ for some
$s,t\in z^{2}M$ and $x\in M$, then
$sz^{2k}=tz^{2k}$, for some $k$.  We may assume
$k$ is even. If $4$ divides $n$ then
$sz^{2k}=s(a_{1}\cdots a_{n})^{k}(a_{n}\cdots
a_{1})^{k}= sa_{1}^{2k}\cdots a_{n}^{2k}$ by
Lemma~\ref{squares}$(ii)$. If $n=4j+2$ for some
$j$ then $sz^{2k}=s(a_{1}a_{2}(a_{1}\cdots
a_{n})(a_{3}\cdots
a_{n}))^{k}=s(a_{1}a_{2}(a_{1}a_{2}a_{n}\cdots
a_{3})a_{3}\cdots
a_{n})^{k}=s(a_{1}a_{2}a_{1}a_{2}a_{3}^{2}\cdots
a_{n}^{2})^{k}=sa_{1}^{2k}\cdots a_{n}^{2k}$, by
Lemma~\ref{ij} and Lemma~\ref{squares}$(ii)$. If
$n=4j+1$ then $z$ is central and thus we get that
$sz^{2k} = s(a_{1}(a_{1}\cdots a_{n})a_{2}\cdots
a_{n})^{k}$  $ = s(a_{1}(a_{1}a_{n}\cdots
a_{2})a_{2}\cdots a_{n})^{k} = sa_{1}^{2k}\cdots
a_{n}^{2k}$, by Lemma~\ref{squares}$(ii)$.
Finally, if $n=4j+3$ then, again by
Lemma~\ref{squares}$(ii)$,
$sz^{2k}=s(a_{1}a_{2}a_{3}(a_{1}\cdots
a_{n})a_{4}\cdots
a_{n})^{k}=s(a_{1}a_{2}a_{3}(a_{3}a_{1}a_{2}a_{4}\cdots
a_{n})a_{n}\cdots
a_{4})^{k}=s(a_{1}a_{2}a_{1}a_{2}a_{3}^{2}\cdots
a_{n}^{2})^{k}= sa_{1}^{2k}\cdots a_{n}^{2k}$.
So, by Lemma~\ref{normalf} and by the previous
comments,  we always get $s=t$. This and a
symmetric argument show that $z^{2}M$ is
cancellative and also that the ideal $z^{2}M$
embeds into $M/\rho$. Hence, again by
Lemma~\ref{canc}, $G= (z^{2})M \langle z^{-2}
\rangle$.

Since $K[G]$ is a PI-algebra and $G$ is the group
of fractions of $M/\rho$ by Lemma~\ref{canc},
$K[M/\rho]$ is a finitely generated PI-algebra.
Let $\overline{M}=M/\rho$.  Since $G$ is a
nilpotent group,
from \cite[Theorem~4.3.3]{bookspringer} and the
comment following it we know that
$K[\overline{M}]$ is noetherian. By
\cite[Theorem~18.1]{book},
$\mathcal{J}(K[\overline{M}])$ is nilpotent.
Therefore, there exists a positive integer $m$
such that $\mathcal{J}(K[M])^m\subseteq I(\rho)$.
By Proposition~\ref{radicalgeneral},
$\mathcal{J}(K[M])^3\subseteq K[z^2M]$. Since
$z^2M$ is cancellative, $I(\rho)\cap K[z^2M]= 0$.
Hence $\mathcal{J}(K[M])$ is nilpotent.

{\it Proof of $(iv)$.}  Suppose that $n\geq 4$ is
odd. We shall see that $a_{1}a_{1}a_{2}z\neq
a_{2}a_{1}a_{1}z$ in $M$. In fact, the only words
beginning with $a_{2}$ that represent the element
$a_{1}a_{1}a_{2}z$ are of the form
$ua_{1}a_{1}a_{2}$, where $u$ is a word beginning
with $a_2$ that represents $z\in M$.

Let $w_0=a_{1}a_{1}a_{2}a_{1}\cdots a_{n}\in \FM_{n}$ and let
$w\in \FM_{n}$ be a word representing the element
$a_{1}a_{1}a_{2}z\in M$. Then there exist $w_1,\dots ,w_r\in
\FM_{n}$ with $w_r=w$ and
$w_{i}=w_{1,i}w_{2,i}w_{3,i}=w'_{1,i}w'_{2,i}w'_{3,i}$ such that
$w_{2,i}$ and $w'_{2,i}$ represent the element $z$ in $M$, for all
$i=0,1,\dots, r$, and $w_{1,j}=w'_{1,j+1}$ and
$w_{3,i}=w'_{3,i+1}$, for all $j=0,\dots ,r-1$. We shall prove, by
induction on $r$, that $w_{1,i}w_{3,i}=a_{1}a_{1}a_{2}$ for all
$i=0,1,\dots ,r$. It is clear that $w_{1,0}=a_{1}a_{1}a_{2}$ and
$w_{3,0}=1$, thus $w_{1,0}w_{3,0}=a_{1}a_{1}a_{2}$. Suppose that
$i\geq 0$ and $w_{1,i}w_{3,i}=a_{1}a_{1}a_{2}$. Then $w_{1,i}\in
\{ 1, a_{1},a_{1}a_{1}, a_{1}a_{1}a_{2}\}$. We shall deal with
four cases separately.

{\em Case $1$:} $w_{1,i}=1$. In this case,
$w_{3,i}=w'_{3,i+1}=a_{1}a_{1}a_{2}$. Since
$w_{i+1}=w_{1,i+1}w_{2,i+1}w_{3,i+1}=w'_{2,i+1}a_{1}a_{1}a_{2}$
and $w_{2,i+1}$ and $w'_{2,i+1}$ represent $z\in
M$, we have that $w_{3,i+1}\in \{
a_{1}a_{1}a_{2}, a_{1}a_{2}\}$. If
$w_{3,i+1}=a_{1}a_{1}a_{2}$, then clearly
$w_{1,i+1}=1$ and
$w_{1,i+1}w_{3,i+1}=a_{1}a_{1}a_{2}$. Suppose
that $w_{3,i+1}=a_{1}a_{2}$. Since the degree in
$a_1$ of $w_{i+1}$ is $3$ and the degree in $a_1$
of $w_{2,i+1}$ is $1$, we have that
$w_{1,i+1}=a_1$. Hence
$w_{1,i+1}w_{3,i+1}=a_{1}a_{1}a_{2}$ in this
case.

{\em Case $2$:} $w_{1,i}=a_{1}$. In this case,
$w_{3,i}=w'_{3,i+1}=a_{1}a_{2}$. Since
$w_{i+1}=w_{1,i+1}w_{2,i+1}w_{3,i+1}=a_{1}w'_{2,i+1}a_{1}a_{2}$,
we have that either $w_{1,i+1}=1$ or $w_{1,i+1}$
begins with $a_1$. If $w_{1,i+1}=1$ then, using
the degree in $a_1$ and that $w_{3,i+1}$ finishes
with $a_{1}a_{2}$, we see that
$w_{3,i+1}=a_{1}a_{1}a_{2}$ and
$w_{1,i+1}w_{3,i+1}=a_{1}a_{1}a_{2}$. Suppose
that $w_{1,i+1}$ begins with $a_1$. Then
$w_{1,i+1}=a_{1}u$ for some $u\in \FM_{n}$. Thus
$uw_{2,i+1}w_{3,i+1}=w'_{2,i+1}a_{1}a_{2}$.  Now
$w_{3,i+1}\in\{ 1, a_2, a_1a_2\}$. If
$w_{3,i+1}\in\{ a_2, a_1a_2\}$, then using the
degree in $a_1$, we have that
$uw_{3,i+1}=a_1a_2$. Suppose that $w_{3,i+1}=1$.
Then $uw_{2,i+1}=w'_{2,i+1}a_1a_2$ and, using the
degree in $a_1$ and in $a_2$, we have that $u\in
\{ a_1a_2, a_2a_1\}$. Since $w_{2,i+1}$ and
$w'_{2,i+1}$ represent $z\in M$,
$f(a_2a_1w_{2,i+1})=1$ and
$f(w'_{2,i+1}a_1a_2)=-1$, where $f$ is the map
defined by~(\ref{f}), thus $u=a_1a_2$. Hence
$w_{1,i+1}w_{3,i+1}=a_{1}a_{1}a_{2}$ in this
case.

{\em Case $3$:} $w_{1,i}=a_{1}a_{1}$. In this
case, $w_{3,i}=w'_{3,i+1}=a_{2}$. Since
$w_{i+1}=w_{1,i+1}w_{2,i+1}w_{3,i+1}=a_{1}a_{1}w'_{2,i+1}a_{2}$
and $w_{2,i+1}$ represents $z\in M$, we have that
$w_{1,i+1}$ begins with $a_1$.  Then
$w_{1,i+1}=a_{1}u$ for some $u\in \FM_{n}$, and
$uw_{2,i+1}w_{3,i+1}=a_1w'_{2,i+1}a_{2}$. Thus,
using the degree in $a_1$ and in $a_2$, we have
$uw_{3,i+1}=a_{1}a_{2}$. Hence
$w_{1,i+1}w_{3,i+1}=a_{1}a_{1}a_{2}$ in this
case.

{\em Case $4$:} $w_{1,i}=a_{1}a_{1}a_2$. In this
case, $w_{3,i}=w'_{3,i+1}=1$. Since
$w_{i+1}=w_{1,i+1}w_{2,i+1}w_{3,i+1}=a_{1}a_{1}a_{2}w'_{2,i+1}$
and $w_{2,i+1}$ represents $z\in M$, we have that
$w_{1,i+1}$ begins with $a_1$.  Then, as in Case
$2$, $w_{1,i+1}=a_{1}u$ for some $u\in \FM_{n}$,
and $uw_{3,i+1}=a_{1}a_{2}$. Hence
$w_{1,i+1}w_{3,i+1}=a_{1}a_{1}a_{2}$ in this
case.

Therefore, we indeed have shown in each of the
four cases that $w_{1,i}w_{3,i}=a_{1}a_{1}a_{2}$,
for all $i=0,1,\dots ,r$. In particular, we have
that $a_{1}a_{1}a_{2}z\neq a_{2}a_{1}a_{1}z$ in
$M$.

Note that if $1\leq i,j\leq n$ are different then
there exists $\sigma\in \Alt_n$ such that
$\sigma(1)=i$ and $\sigma(2)=j$. Therefore
$$a_{i}a_{i}a_{j}z\neq a_{j}a_{i}a_{i}z,$$ for
all $i\neq j$, in $M$.

Since $n$ is odd,  $z$ is central in $M$ and, by
Lemma~\ref{squares},
$(a_ia_ia_j-a_ja_ia_i)z^{2}=0$. Therefore
$(a_ia_ia_j-a_ja_ia_i)z\in
\mathcal{B}(K[M])\setminus \{0\}$, for all $i\neq
j$ and for any field $K$.

Let $\overline{\rho}=\rho\cap (zM\times zM)$. So
$I(\overline{\rho})=\lin_K \{ s-t\mid s,t\in zM \mbox{ and
}\exists i\geq 0, \, sz^{i}=tz^{i}\}$. Since $z^{2}M$ is
cancellative, it follows that $I(\overline{\rho})^2=0$.

Suppose that $K$ has characteristic $\neq 2$. We
have that $\mathcal{J}(K[G])=0$. Since
$\mathcal{J}(K[\overline{M}])$ is nilpotent and
$G$ is a central localization of $\overline{M}$,
we get
$\mathcal{J}(K[\overline{M}])=\mathcal{B}(K[\overline{M}])\subseteq
\mathcal{J}(K[G])$. Hence
$\mathcal{J}(K[\overline{M}])=0$. Then, by
Corollary~\ref{radical},
$$\mathcal{J}(K[M])=\mathcal{B}(K[M])=I(\overline{\rho}).$$
Thus $\mathcal{J}(K[M])^2=0$.

Assume that $s,t\in M$ are such that $(s,t)\in \rho$. Because
$z^{2}M$ is cancellative, we know that $z^{2}s=z^{2}t$. From the
proof of Lemma~\ref{normalf} it follows that $s$ can be rewritten
in finitely many steps, each step being a rewriting of one of the
following types: $a_{i}a_{j}a_{k}\mapsto a_{j}a_{k}a_{i}$,
$a_{i}a_{j}^{2}\mapsto a_{j}^{2}a_{i}$ and
$a_{i}a_{j}a_{i}a_{j}\mapsto a_{j}a_{i}a_{j}a_{i}$, so that
$z^{2}s$ is in one of the forms given in this lemma. The same
applies to the element $z^{2}t$. Since $z^{2}s=z^{2}t$, a degree
argument shows then that these new forms of $z^{2}s$ and $z^{2}t$
are equal. Hence the forms of $s$ and $t$ are equal. It follows
that $s-t\in K[M]YK[M]$, where
$$Y=\bigcup_{\stackrel{1\leq i,j,k\leq n}{|\{
i,j,k \} |=3}} \{ a_ia_ja_k-a_ja_ka_i,\; a_i^2a_j-a_ja_i^2,\;
(a_ia_j)^2-(a_ja_i)^2\}.$$ This implies that $Y$ generates
$I(\rho)$ as a  two-sided ideal. Now, if $s'z,t'z\in zM$ are
$\rho$-related, then also $(s',t')\in \rho$, so by the previous
$s'z-t'z\in K[M]YzK[M]$ because $z$ is central. In particular,
$I(\overline{\rho})=\mathcal{J}(K[M])$ is a finitely generated
ideal.

Suppose that $K$ has characteristic $2$. By
Proposition~\ref{radicalgeneral}, $\mathcal{J}(K[M])\subseteq
K[zM]$. Thus $\mathcal{J}(K[M])=\mathcal{J}(K[zM])$. Note that
$zM/\overline{\rho}$ is a cancellative semigroup and $G$ is its
group of fractions. Furthermore,
$K[zM/\overline{\rho}]=K[zM]/I(\overline{\rho})$. By $(iii)$, we
have that $\mathcal{J}(K[M])$ is nilpotent. Hence
$$\mathcal{J}(K[M])/I(\overline{\rho})=\mathcal{B}(K[zM])/I(\overline{\rho})
=\mathcal{B}(K[zM/\overline{\rho}])=\mathcal{B}(K[G])\cap
K[zM/\overline{\rho}],$$  see \cite[Corollary~11.5]{book}. Let
$\pi\colon zM\longrightarrow G/C$ be the composition of the
natural maps
$$zM\hookrightarrow M\longrightarrow
G\longrightarrow G/C.$$ Let $\eta$ be the congruence defined on
$M$ by $s\eta t$ if and only if either $s=t$ or $s,t\in zM$ and
$\pi(s)=\pi(t)$. Since $\mathcal{B}(K[G])=\omega (K[C])K[G]$, it
follows that $\mathcal{J}(K[M])=I(\eta)$.  In particular,
$z(a_{i}a_{j}-a_{j}a_{i})\in \mathcal{J}(K[M])$ for all $i,j$. Let
$Q$ be the ideal of $K[M]$ generated by all such elements. Then
the set of all elements of the form $za_{1}^{i_{1}}\cdots
a_{n}^{i_{n}}$, for nonnegative integers $i_{1},\ldots ,i_{n}$,
forms a basis of the algebra $K[zM]/Q$. Therefore it embeds into
the algebra $K[G/C]$, which is a commutative domain. It follows
that ${\mathcal J}(K[M])=Q$ and hence it is finitely generated.

{\it Proof of $(v)$.} Let $K$ be any field. Suppose
that $n\geq 6$ is and $n$ is even.  We shall prove
that $\mathcal{J}(K[M])\subseteq K[T]$ (where $T$ is
as in Lemma~\ref{even2}). Suppose that
$\mathcal{J}(K[M])\not\subseteq K[T]$. Let
$\alpha\in \mathcal{J}(K[M])\setminus K[T]$ with
$|\supp(\alpha)|=m$. Let $\supp(\alpha)=\{ s_1,\dots
,s_m\}$. Suppose that $s_1\notin T$. Then, by
Lemma~\ref{even2}, there exist $i,j$ such that
$s_1a_j\notin Mz$ and $a_is_1\notin zM$. Hence
$a_is_1a_j\notin zM\cup Mz$
 and $a_is_1a_j\in\supp
(a_{i}\alpha a_{j})$. But this is in contradiction
with Lemma~\ref{radicalgeneral}. Hence
$\mathcal{J}(K[M])\subseteq K[T]$.

We now prove that $K[T]\cap I(\rho)=0$, i.e., $T$
is cancellative. Let $s,t\in T$ be such that
$s\rho t$. By Lemma~\ref{normalf2} $(i)$ and
$(ii)$, there exist $s',t'\in M$, three different
integers $1\leq i,j,k\leq n$ and three different
integers $1\leq l,p,q\leq n$ such that
$s=a_ia_ja_kzs'$ and $t=a_la_pa_qzt'$. Suppose
that $l\notin\{ i,j,k\}$. Since $s\rho t$, $s'\in
Ma_lM$. By Lemma~\ref{ij},
$(a_ia_ja_kz)a_l=a_iza_ja_ka_l=a_i(a_ka_ja_l)z$.
By Lemma~\ref{normalf2} $(i)$ and $(ii)$, we may
assume that $\{i,j,k\}=\{l,p,q\}$. Suppose that
$s=a_ia_ja_kzs'$ and $t=a_ia_ja_kzt'$. In this
case, by Lemma~\ref{normalf2}, $s=t$. Thus, by
Lemma~\ref{ij}, we may assume that
$s=a_ia_ja_kzs'$ and $t=a_ja_ia_kzt'$. Suppose
that $s'\in \bigcup_{1\leq r\leq n}\langle
a_r\rangle$. In this case, since $s$ and $t$ have
the same degrees with respect to every $a_l$,
there exists $1\leq r\leq n$, such that $s'=t'\in
\langle a_r\rangle$, but $a_ia_j\neq a_ja_i$ in
$G$, a contradiction. Therefore, $s'\notin
\bigcup_{1\leq r\leq n}\langle a_r\rangle$. Hence
there exists different $1\leq r,p\leq n$ and
$s_1,s_2\in M$ such that $s'=s_1a_ra_ps_2$. By
Lemma~\ref{normalf2}, there exists a permutation
$\sigma$ of $\{ i,j,k\}$, such that
$s=s_1a_{\sigma(i)}a_{\sigma(j)}a_{\sigma(k)}za_ra_ps_2$.
By Lemma~\ref{ij}, we may assume that
$\sigma(k)\notin\{ r,p\}$. Now, by
Lemma~\ref{ij}, we have
\begin{eqnarray}    \label{fakt}
a_{\sigma(i)}a_{\sigma(j)}a_{\sigma(k)}za_ra_p
&=&
za_{\sigma(j)}a_{\sigma(i)}a_{\sigma(k)}a_ra_p
=a_{\sigma(j)}a_{\sigma(i)}za_{\sigma(k)}a_ra_p\\
&=& a_{\sigma(j)}a_{\sigma(i)}a_{\sigma(k)}a_pa_rz
=a_{\sigma(j)}a_{\sigma(i)}a_{\sigma(k)}za_pa_r. \nonumber
\end{eqnarray}
Hence
$s=s_1a_{\sigma(j)}a_{\sigma(i)}a_{\sigma(k)}za_pa_rs_2$.
By Lemma~\ref{normalf2},
$s=a_ja_ia_kzs_1a_pa_rs_2$, and since
$t=a_ja_ia_kzt'$, $s=t$.

Therefore $K[T]\cap I(\rho)=0$. Hence, if $K$ has
characteristic $\neq 2$, then since
$\mathcal{J}(K[M])\subseteq K[T]$ and
$\mathcal{J}(K[M])\subseteq I(\rho)$ by
Corollary~\ref{radical}, it follows that
$\mathcal{J}(K[M])= 0$.

On the other hand, if $K$ has characteristic $2$
then $\mathcal{J}(K[M])\subseteq K[T]$ implies that
$\mathcal{J}(K[M])=\mathcal{J}(K[T])$. Since $T$ is
a cancellative ideal in $M$ and $G$ is its group of
fractions, in fact $G=T\langle z^{2}\rangle^{-1}$,
by $(iii)$, we have that
$$\mathcal{J}(K[M])=\mathcal{B}(K[M])=\mathcal{B}(K[T])=\mathcal{B}(K[G])\cap
K[T],$$  by \cite[Corollary~11.5]{book}. Let $\pi\colon
T\longrightarrow G/C$ be the composition of the natural maps
$$T\hookrightarrow G\longrightarrow G/C.$$ Let $\eta$ be the
congruence defined on $M$ by $s\eta t$ if and only if either $s=t$
or $s,t\in T$ and $\pi(s)=\pi(t)$. Since $\mathcal{B}(K[G])=\omega
(K[C])K[G]$, it follows that $\mathcal{J}(K[M])=I(\eta)$. Note
that $a_1a_2a_3z+a_2a_1a_3z\in K[T]\cap \mathcal{B}(K[G])$ and it
is nonzero.

Let $P$ be the ideal of $K[M]$ generated by all elements of
the form $a_{j}a_{i}a_{k}z-a_{i}a_{j}a_{k}z$, with different
$i,j,k$. By (\ref{fakt}), $a_{i}a_{j}a_{k}za_{p}a_{q}$ and
$a_{j}a_{i}a_{k}za_{q}a_{p}$ are $\rho$-related, for all $p\neq
q$. By Lemma~\ref{normalf2}, $a_{j}a_{i}a_{k}z$ becomes central in
$K[M]/P$. As at the beginning of the proof of cancellativity of
$T$, if $s\in T$ then $s=a_{i}a_{j}a_{k}zw$ or
$s=a_{j}a_{i}a_{k}zw$, for some $w\in M$, where $i,j,k$ are
minimal indices such that the degree of $s$ in each of
$a_{i},a_{j}$ and $a_{k}$ is positive. Therefore
$s-a_{i}a_{j}a_{k}zw\in P$. Furthermore, we may also assume that
$w$ is of the form $a_{q}^{i_{q}}\cdots a_{n}^{i_{n}}$ for some
$q\geq \min \{ i,j,k\}$ and nonnegative integers $i_{q},\ldots
,i_{n}$. If $t\in T$ is such that $(s,t)\in \eta$, then $s$ and
$t$ have the same degree in each generator and hence $s-t\in P$.
It follows that $\mathcal{J}(K[M])=I(\eta)=P$ and hence it is
finitely generated.

This finishes the proof of
Theorem~\ref{main-alt}.

 \vspace{30pt}
 \noindent \begin{tabular}{llllllll}
 F. Ced\'o && E. Jespers  \\
 Departament de Matem\`atiques &&  Department of Mathematics \\
 Universitat Aut\`onoma de Barcelona &&  Vrije Universiteit Brussel  \\
08193 Bellaterra (Barcelona), Spain    && Pleinlaan 2, 1050
Brussel, Belgium \\ cedo@mat.uab.cat &&
efjesper@vub.ac.be\\
   &&   \\
J. Okni\'nski &&  \\ Institute of Mathematics &&
\\ Warsaw University&& \\ Banacha 2&& \\ 02-097
Warsaw, Poland &&\\ okninski@mimuw.edu.pl
\end{tabular}

\end{document}